\documentclass[11pt,reqno,a4paper]{amsart}

\usepackage[plain]{fullpage}

\usepackage{calrsfs}
\usepackage[OT2,T1]{fontenc}
\usepackage[english]{babel}
\usepackage{newtxtext}
\usepackage{newtxmath}
\usepackage{tikz}
\usetikzlibrary{patterns}
\usetikzlibrary{datavisualization.polar}

\usepackage{etoolbox}
\patchcmd{\section}{\scshape}{\bfseries}{}{}
\makeatletter
\renewcommand{\@secnumfont}{\bfseries}
\makeatother

\usepackage{enumitem}

\theoremstyle{plain}
\newtheorem{theorem}{Theorem}[section]
\newtheorem{proposition}[theorem]{Proposition}
\newtheorem{lemma}[theorem]{Lemma}
\newtheorem{corollary}[theorem]{Corollary}
\theoremstyle{definition}

\counterwithin{equation}{section}

\setenumerate[0]{label=\textup{(\alph*)}}

\usepackage[unicode=true,pdfusetitle, bookmarks=true,bookmarksnumbered=false,bookmarksopen=false, breaklinks=true,pdfborder={0 0 0},pdfborderstyle={},backref=false,colorlinks=true]{hyperref}
\definecolor{myblue}{rgb}{0.09,0.32,0.44} 
\hypersetup{pdfborder={0 0 0},pdfborderstyle={},colorlinks=true,linkcolor=myblue,citecolor=myblue,urlcolor=blue}
\usepackage{mathtools}
\usepackage{mathrsfs}

\newcommand{\Z}{\mathbf{Z}}
\newcommand{\R}{\mathbf{R}}
\newcommand{\Q}{\mathbf{Q}}
\newcommand{\FF}{\mathbf{F}}
\newcommand{\Cc}{\mathbf{C}}
\newcommand{\F}{\mathscr{F}}
\newcommand{\mr}{\mathrm}
\renewcommand{\geq}{\geqslant}
\renewcommand{\leq}{\leqslant}
\DeclareMathOperator{\Sq}{Sq}
\DeclareMathOperator{\erf}{erf}
\newcommand{\ii}{\mathrm{i}}
\newcommand\dif{\mathop{}\!\mathrm{d}} 

\DeclarePairedDelimiter\abs{\lvert}{\rvert}%
\DeclarePairedDelimiter\norm{\lVert}{\rVert}%

\makeatletter
\let\oldabs\abs
\def\abs{\@ifstar{\oldabs}{\oldabs*}}
\let\oldnorm\norm
\def\norm{\@ifstar{\oldnorm}{\oldnorm*}}
\makeatother

\begin{document}

\date{\today}
\title{Counting (skew-)reciprocal Littlewood polynomials \\ with square discriminant}
\author[D.~Hokken]{David Hokken}
\address{\normalfont Mathematisch Instituut, Universiteit Utrecht, Postbus 80.010, 3508 TA Utrecht, Nederland}
\email{d.p.t.hokken@uu.nl}

\subjclass[2020]{Primary: 11C08, 11R32, 11R09, 05A16. Secondary: 11P21.} 
\keywords{\normalfont Littlewood polynomials, square discriminant, Galois theory, asymptotic enumeration, lattice points}

\begin{abstract} \noindent
A \emph{Littlewood polynomial} is a univariate polynomial all of whose coefficients lie in $\{ \pm 1\}$. We establish the leading term asymptotics of the number of reciprocal or skew-reciprocal Littlewood polynomials with square discriminant. This relates to a bounded-height analogue of the Van der Waerden conjecture on Galois groups of random polynomials. As a byproduct, we establish the asymptotics of certain Gaussian-weighted counts of Pythagorean triples.
\end{abstract}

\maketitle

\section{Introduction}

\subsection*{Background and main result}
Let $f$ be a monic polynomial of degree $n$ with integer coefficients that are at most $H$ in absolute value. In 1934, Van der Waerden \cite{VdW1} presented an elementary proof that $f$ is almost surely \emph{ohne Affekt}: the Galois group $G_f$ of $f$ over $\Q$ is the symmetric group $S_n$ with probability tending to $1$ as $H$ goes to infinity. Two years later, he posed a conjecture \cite[p.~139]{VdW2} on the probability that $f$ does not have maximal Galois group, which states
\begin{equation} \label{eq:VdW}
\mr{Prob}(G_f \neq S_n) \sim \mr{Prob}(f \text{ is reducible})
\end{equation}
as $H$ goes to infinity. Last year, Bhargava \cite{Bhargava} established the breakthrough result  
\begin{equation} \label{eq:Bhargava}
\mr{Prob}(G_f \neq S_n) \sim \mr{Prob}(f \text{ is reducible}) + \mr{Prob}(G_f = A_n) \asymp H^{-1}
\end{equation}
where $A_n$ denotes the alternating group on $n$ letters. This is a weak form of the Van der Waerden conjecture.
Since $f$ is reducible with probability $\asymp H^{-1}$ if $n > 2$ (see \cite{VdW1, Chela}), the remaining task to obtain \eqref{eq:VdW} consists of showing that $\mr{Prob}(G_f = A_n) = o(H^{-1})$; Bary-Soroker, Ben-Porath and Matei \cite{BBM} conjecture the much stronger bound ${\mr{Prob}(G_f = A_n) = O(H^{-n/2+\epsilon})}$ when $n \geq 4$.

The \emph{height} $H$ of the polynomial $f$ in the above setup tends to infinity, whereas the degree $n$ stays fixed. This approach to random polynomials is called the \emph{large box model}. In the \emph{restricted coefficient model}, the height $H$ --- or any specific set $\mathcal{N}$ of coefficients of $f$ --- is fixed, and it is the degree that tends to infinity. Recent years have seen a surge of interest in questions about Galois groups in this setting as well \cite{BKK, BK, Borst_et_al, BV, Konyagin, OW}. For example, if $f$ is sampled uniformly at random among the monic degree-$n$ polynomials all of whose coefficients lie in a set $\mathcal{N}$ of at least $35$ consecutive integers and whose constant coefficient is nonzero, Bary-Soroker, Koukoulopoulos and Kozma \cite{BKK} show that $G_f$ is $S_n$ or the alternating group $A_n$ with probability tending to $1$. Conditionally on the Riemann Hypothesis for a family of Dedekind zeta functions, Breuillard and Varj\'{u} \cite{BV} show a similar result for more general distributions of the coefficients of $f$.
The methods from the large box model yield a dependence of the implicit constant in \eqref{eq:Bhargava} on the degree $n$, meaning that they do not apply in the restricted coefficient model. Instead, probabilistic methods and/or finite group theory are used to establish high transitivity of $G_f$ from the reductions of $f$ modulo various primes.
This leaves only $A_n$ and $S_n$ as possible Galois groups, but as these are respectively $(n-2)$- and $n$-transitive, it is hard to distinguish them based on this property. In other words, the alternating group has a special role in the restricted coefficient model as well. Generally, it is believed that $A_n$ should occur with probability tending to $0$ as $n$ tends to infinity \cite{BK}.

Recall that the \emph{discriminant} of the monic polynomial $f$ is the integer $\Delta(f) := \prod_{i < j} (\alpha_i-\alpha_j)^2$, where $\alpha_1, \ldots, \alpha_n$ denote the complex roots of $f$ (with multiplicities). A property that distinguishes $A_n$ from $S_n$ as Galois group $G_f$ of a separable polynomial $f$ is the following: $G_f$ is contained in $A_n$ if and only if the discriminant $\Delta(f)$ of $f$ is a (necessarily nonzero) square. This implies the purely algebraic fact formulated in probabilistic terms as $\mr{Prob}(G_f = A_n) \leq \mr{Prob}(\Delta(f) = \square \neq 0)$, independent of the underlying distribution on the set of polynomials.

This paper studies the probability that the discriminant of the monic polynomial $f$ is a square when the coefficients of $f$ are independently and uniformly selected from $\{\pm 1\}$. Such polynomials are often called \emph{Littlewood polynomials}. These are extremal examples of polynomials with restricted coefficients: all Littlewood polynomials in degree $n$ coincide over $\FF_2$, whereas they form a sparse (that is, exponentially small in $n$) subset of the degree-$n$ monic polynomials in $\FF_p[X]$ for any prime $p > 2$. Furthermore, since they are of height $1$, the results mentioned in the first paragraph cannot be made effective in any way. The state-of-the-art result concerning the Galois theory of random Littlewood polynomials is that at least a fraction of $0.00068$ of the Littlewood polynomials of degree $n$, with $n \geq 10^{10^{4.9}}$, is irreducible (see \cite[Theorem~3.5]{BKK}).

Following Littlewood \cite{Littlewood}, denote the collection of Littlewood polynomials of degree $n$ by $\F_n$; let $\Sq_n \subset \F_n$ consist of those with square discriminant. Furthermore, call $f$ \emph{reciprocal} if $f(X) = X^{n} f(X^{-1})$ and \emph{skew-reciprocal} if $f(X) = (-1)^{n(n-1)/2}X^n f(-X^{-1})$ (the latter appear e.g. in \cite{Odlyzko, Erdelyi} in connection to questions about the flatness of Littlewood polynomials on the unit circle). Denote by $R_n, S_n \subset \F_n$ the sets of Littlewood polynomials of degree $n$ that have square discriminant and are reciprocal, respectively skew-reciprocal. Our main result concerns the size of $R_n$ and $S_n$ as $n$ tends to infinity.
\begin{theorem} \label{thm:mainthm} 
The sets $R_{8n}$, $S_{8n}$, $R_{8n-2}$, and $S_{8n-2}$ are all of size $\asymp 16^{n} \log{n}/\sqrt{n}$. More precisely:
\begin{enumerate}
    \setlength{\itemsep}{5pt}
    \item $\displaystyle \lim_{n \to \infty}{\frac{\abs{R_{8n}}}{16^n \log{n}/\sqrt{n}}} = \frac{\Gamma(\frac{1}{4})^2}{4\sqrt{2} \pi^3} = 0.0749\ldots,$ \\[5pt]
    $\displaystyle \lim_{n \to \infty}{\frac{\abs{S_{8n}}}{16^n \log{n}/\sqrt{n}}} = \frac{1}{2\pi^{3/2}} = 0.0897\ldots;$
    \item $\displaystyle \abs{R_{8n-2}} \sim \frac{1}{2} \abs{R_{8n}}$ and $\displaystyle \abs{S_{8n-2}} \sim \frac{1}{2} \abs{S_{8n}}$.
\end{enumerate}
\end{theorem}
The limits in Theorem \ref{thm:mainthm} are approached extremely slowly. For example, when $n = 10^{11}$, the fraction $\abs{R_{8n}}/(16^n \log{n}/\sqrt{n})$ is $0.099\ldots$. This is (at least in part) due to large contributions of order $\asymp 16^n/\sqrt{n}$ to $\abs{R_{8n}}$ and $\abs{S_{8n}}$ coming from error terms in lattice point counts that we use.

As observed in \cite[\S 4]{BK}, any $f \in \F_{2n}$ of even degree is separable, because $f$ coincides modulo $2$ with the separable polynomial $(X^{2n+1}-1)/(X-1)$. Furthermore, the roots of a reciprocal polynomial $f$ come in pairs $\{\alpha, \alpha^{-1}\}$; if $f$ is skew-reciprocal, they come in pairs $\{\alpha, -\alpha^{-1}\}$. The separability of $f$ implies that $\alpha$ and $\pm \alpha^{-1}$ are distinct. As a result, the Galois group of (skew-)reciprocal $f$ is contained in the \emph{permutational wreath product} $C_2 \wr S_{n/2}$, see \cite{VV}. Recall that the permutational wreath product of two groups $G$ and $H \leq S_n$, denoted $G \wr H$, is the semidirect product $G^n \rtimes H$ where $H$ acts on the $n$ copies of $G$ by permuting the coordinates. Theorem \ref{thm:mainthm} therefore leads to the following corollary.

\begin{corollary} \label{cor:maincor}
Let $f$ be sampled uniformly at random among the (skew-)reciprocal Littlewood polynomials of degree $n \equiv 0, 6 \bmod{8}$. As $n \to \infty$, we have
\begin{equation*}
\mr{Prob}(\Delta(f) = \square \neq 0) = \mr{Prob}(G_f \leq (C_2 \wr S_{n/2}) \cap A_{n}) \asymp \frac{\log{n}}{\sqrt{n}},
\end{equation*}
the implicit constants being as specified in Theorem~\ref{thm:mainthm}.
\end{corollary}

The set $\Sq_n$ is empty whenever $n \equiv 2, 4 \bmod{8}$, which is the reason to leave out these degrees in the above statements. In \S \ref{sec:otherdeg}, we expound the proof sketch for this fact provided in \cite[\S 4]{BK}. In the same section we also make some remarks on the case of odd $n$.

Reciprocals and skew-reciprocals are decomposable: a polynomial $f$ is reciprocal if it is of the form $f(X) = X^{n/2} g(X + X^{-1})$, and skew-reciprocal if it is of the form $f(X) = X^{n/2} g(X - X^{-1})$ for some polynomial $g$. The group $(C_2 \wr S_{n/2}) \cap A_{n}$ is much smaller than $A_n$ --- of index $1\cdot 3 \cdot 5 \cdot \ldots \cdot (n-1)$ to be precise --- and the sizes of $R_n$ and $S_n$ compared to $\abs{\F_n} = 2^n$ decrease exponentially in $n$. Nevertheless, back in the large box model, the best known bound on the probability that the discriminant of $f$ is a square also come from decomposable polynomials with the very same Galois group: Bary-Soroker, Ben-Porath and Matei \cite[Theorem~1.3]{BBM} show for all even $n \geq 6$ that
\begin{equation*}
\mr{Prob}(\Delta(f) = \square) \gg H^{-(n+1)/2}
\end{equation*}
as $H$ tends to infinity by applying an explicit version of Hilbert's irreducibility theorem to polynomials of the form $f(X) = g(X^2)$.
No Littlewood polynomial of the form $f(X) = g(X^2)$ exists, and it appears that (skew-)reciprocal polynomials are `the next best thing'. 

\subsection*{Outline} In the setting of Littlewood polynomials, reducing modulo primes or applying probabilistic methods seems difficult. Instead, the proofs in this paper combine counting arguments to derive explicit combinatorial expressions for the objects of study with lattice point counts in certain geometric regions and asymptotics of binomial coefficients.

In \S \ref{sec:counting} we derive combinatorial expressions for $\abs{R_{8n}}$ and the three other sets under consideration, see Proposition \ref{prop:nrofsqdisc} and Proposition \ref{prop:nrofsqdiscother}. In each case, we obtain a sum that extends over certain tuples related to Pythagorean triples; these come from a square discriminant criterion for (skew-)reciprocal polynomials given in Lemma \ref{lem:discrit}. This criterion can in theory be used to find similar expressions when Littlewood polynomials are replaced by polynomials with coefficients in any fixed set $\mathcal{N}$. Auxiliary results to study the asymptotics of these combinatorial expressions, as well as an analysis of the Pythagorean triples, are contained in \S \ref{sec:aux}. The latter essentially boils down to counting lattice points with parity and coprimality conditions inside elliptic (for the reciprocals) or parabolic (for the skew-reciprocals) hyperboloids. These results are then combined in \S \ref{sec:reciprocals} and \S \ref{sec:skewreciprocals}, where the lattice point regions are split into three suitably chosen parts. This makes it possible to evaluate the combinatorial expressions from \S \ref{sec:counting} asymptotically by using integral estimates. The proof of Theorem \ref{thm:mainthm} is finally given in \S \ref{sec:mainthm}. In \S \ref{sec:pyth}, we discuss implications of our results on certain Gaussian-weighted counts of Pythagorean triples, see Theorem~\ref{thm:pyth} below. We end with some observations about the set $\Sq_n$ in the case $n \not \equiv 0, 6 \bmod{8}$ in \S \ref{sec:otherdeg}. 

\subsection*{Counting Pythagorean triples}
As touched upon in the outline, it turns out that there is a relation between counting (skew-)reciprocal Littlewood polynomials with square discriminant and counting Pythagorean triples. This relation is explained in \S \ref{sec:counting}, and provides combinatorial expressions for the main objects in this paper. As a byproduct of the results in this paper, we establish the following theorem.
\begin{theorem}
\label{thm:pyth}
Suppose $A_0, A_1, \ldots$ are random variables taking the values $\pm 1$ with probability $1/2$ each, and set $X = 2\sum_{i=1}^{2n} A_{2i-1}$ and $Y = A_0 + 2\sum_{i=1}^{2n} A_{2i}$. As $n$ goes to infinity, we have
\begin{equation} \label{eq:probgauss}
\mr{Prob}(Y^2-X^2 = \square) \sim \frac{\Gamma(\frac{1}{4})^2}{4\sqrt{2}\pi^3} \frac{\log{n}}{\sqrt{n}} \qquad \text{and} \qquad \mr{Prob}(Y^2+X^2 = \square) \sim \frac{1}{2\pi^{3/2}} \frac{\log{n}}{\sqrt{n}}.
\end{equation}
\end{theorem}
The proof of Theorem~\ref{thm:pyth} is given in \S \ref{sec:pyth}.
Note that $X$ and $Y$ are simple random walks with step size $2$ (with $Y$ starting with a single step of size $1$). In the limit as $n$ tends to infinity, the relevant local central limit theorem \cite[Theorem~3.1.2]{Durrett} implies that, locally and close to the origin, the random variables $X$ and $Y$ tend to Gaussians with standard deviations on the order of $\sqrt{n}$. In fact, we can think of Theorem~\ref{thm:pyth} as the asymptotics of certain Gaussian-weighted counts of Pythagorean triples, see \eqref{eq:Gauss}. 

In the unweighted case, the following are results due to Sierpiński \cite[Chapter~5,~Eq.~(7)]{Sierpinski} and Benito and Varona \cite[Corollary~2]{BVPyth}, respectively: if $X$ and $Y$ are independent random variables taking values in $[-\sqrt{n}, \sqrt{n}] \cap \Z$ uniformly, then
\begin{equation} \label{eq:probunif}
\mr{Prob}(Y^2-X^2 = \square) \sim \frac{1}{2\pi} \frac{\log{n}}{\sqrt{n}} \qquad \text{and} \qquad \mr{Prob}(Y^2+X^2 = \square) \sim \frac{2 \log(1+\sqrt{2})}{\pi^2} \frac{\log{n}}{\sqrt{n}}.
\end{equation}
The constants $1/(2\pi)$ and {$2\log(1+\sqrt{2})/\pi^2$} in \eqref{eq:probunif} should be divided by $2$ in case we assume, as in Theorem~\ref{thm:pyth}, that $Y$ is odd and $X$ is even. Notice that the asymptotic sizes in \eqref{eq:probgauss} and \eqref{eq:probunif} are all on the order of $\log{n}/\sqrt{n}$. It would be interesting to understand if it is possible to pass more directly from the classical unweighted results, which we use here as well (see Proposition~\ref{prop:Bnsize} and Proposition~\ref{prop:Dnsize}), to our result for the Gaussian-weighted case.

\subsection*{Notation}

The expression $f \ll g$ as well as $g \gg f$ and $f = O(g)$ all mean there exists a positive constant $C$ such that $\abs{f(n)} \leq C \abs{g(n)}$ for all sufficiently large values of $n$ (all asymptotics in this paper will be in $n$). The notation $f \asymp g$ is shorthand for $g \ll f \ll g$. The functions $f$ and $g$ are said to be \emph{asymptotically equal}, denoted $f \sim g$, if the fraction $f(n)/g(n)$ tends to $1$ as $n$ tends to infinity. In particular, $f \sim g$ implies $f \asymp g$. Lastly, the notation $f = o(g)$ is used when the fraction $f(n)/g(n)$ tends to $0$ as $n$ tends to infinity.

We write $i$ for an index and $\ii \in \Cc$ for the imaginary unit and adopt the convention that $\binom{n}{k} = 0$ if $k > n$.

\subsection*{Acknowledgements}
Many thanks to Gunther Cornelissen, Mar Curc\'{o} Iranzo, and Berend Ringeling for helpful conversations and feedback on earlier versions of this manuscript. The author thanks the two anonymous reviewers for their helpful comments and interesting questions that greatly improved the paper. This publication is part of the project \emph{Littlewood polynomials with square discriminant} (OCENW.M20.233), financed by the Dutch Research Council (NWO).

\section{A counting proposition}
\label{sec:counting}

In this section, we prove the following expression for $\abs{R_{8n}}$ in terms of binomial coefficients.
\begin{proposition} \label{prop:nrofsqdisc}
The number of reciprocal Littlewood polynomials of degree $8n$ with (nonvanishing) square discriminant equals
\begin{equation} \label{eq:nrofsqdisc}
\abs{R_{8n}} = 2^{2n} \binom{2n}{n} + 2 \sum \binom{2n}{n+\frac{1}{2}krs} \binom{2n}{n+\frac{1}{4}(kr^2+ks^2 + (-1)^{\frac{k+1}{2}})}
\end{equation}
where the sum extends over all tuples $(k,r,s)$ such that $k>0$ is odd and $r>s>0$ are coprime and of opposite parity (i.e., $r$ is odd if and only if $s$ is even).
\end{proposition}
Similar expressions for $\abs{R_{8n-2}}$, $\abs{S_{8n}}$ and $\abs{S_{8n-2}}$ are given in Proposition \ref{prop:nrofsqdiscother}. The first term in \eqref{eq:nrofsqdisc} is $\asymp 16^n/\sqrt{n}$ as a consequence of the well-known asymptotic expression ${\binom{2n}{n} \sim 4^n/\sqrt{\pi n}}$ for the central binomial coefficient \cite[\S 5.4]{Spencer}. Theorem \ref{thm:mainthm} claims that this falls short by a factor logarithmic in $n$ of the true growth rate.

The proof of Proposition \ref{prop:nrofsqdisc} is based on the following square discriminant criterion.
 
\begin{lemma}
\label{lem:discrit}
Let $f \in \Q[X]$ be a separable polynomial of degree $2n$. Suppose $f$ is reciprocal. Then the discriminant of $f$ is a square if and only if $(-1)^n f(1)f(-1)$ is a square. Similarly, if $f$ is skew-reciprocal, then its discriminant is a square if and only if the integer $f(\ii)f(-\ii)$ is a square.
\end{lemma}
\begin{proof}
In the case of reciprocal polynomials, this criterion is well-known and recorded in the literature in several places, see e.g. \cite[p.~85]{Dubickas}. With a similar proof, here we show the criterion for skew-reciprocals.

Write $a_n$ for the leading coefficient of $f$. If $f$ is not monic, then $\Delta(f) = a_n^{2n-2} \Delta(f/a_n)$. Since $a_n^{2n-2}$ is a square, we may assume without loss of generality that $f$ is in fact monic. Since $f$ is separable, it has $2n$ distinct roots. These come in pairs $\alpha_i, \alpha_{n+i} = -\alpha_i^{-1}$ for $i = 1, \ldots, n$. Hence
\begin{equation*}
\Delta(f) = \prod_{1 \leq i < j \leq n} \left((\alpha_i - \alpha_j)(\alpha_i+\alpha_j^{-1})(-\alpha_i^{-1}+\alpha_j^{-1})(-\alpha_i^{-1}-\alpha_j) \right)^2 \prod_{1 \leq j \leq n} (\alpha_j + \alpha_j^{-1})^2.
\end{equation*}
The first of the two products above is the square of an integer, since
\begin{equation*}
\prod_{1 \leq i < j \leq n} (\alpha_i - \alpha_j)(\alpha_i+\alpha_j^{-1})(-\alpha_i^{-1}+\alpha_j^{-1})(-\alpha_i^{-1}-\alpha_j) = \prod_{1 \leq i < j \leq n} -(\alpha_i-\alpha_i^{-1} - \alpha_j+\alpha_j^{-1})^2
\end{equation*}
is a symmetric expression in the roots of $f$.
The other product can be expanded as
\begin{equation*}
\prod_{1 \leq j \leq n} (\alpha_j + \alpha_j^{-1})^2 = \prod_{1 \leq j \leq n} (\ii+\alpha_j)(\ii+\alpha_j^{-1})(\ii - \alpha_j)(\ii-\alpha_j^{-1}) =  f(\ii)f(-\ii)
\end{equation*}
as claimed.
\end{proof}

To count (skew-)reciprocal polynomials with square discriminant, we recall that any polynomial $f$ can be written as the sum $f(X) = f_{\mr{e}}(X^2) + Xf_{\mr{o}}(X^2)$ of its even and odd parts. Therefore
\begin{equation} \label{eq:countrec}
f(1)f(-1) = (f_{\mr{e}}(1) + f_{\mr{o}}(1))(f_{\mr{e}}(1) - f_{\mr{o}}(1)) = f_{\mr{e}}(1)^2 - f_{\mr{o}}(1)^2
\end{equation}
and
\begin{equation} \label{eq:countskew}
f(\ii)f(-\ii) =  (f_{\mr{e}}(\ii^2) + \ii f_{\mr{o}}(\ii^2))(f_{\mr{e}}((-\ii)^2)-\ii f_{\mr{o}}((-\ii)^2)) = f_{\mr{e}}(-1)^2 + f_{\mr{o}}(-1)^2.
\end{equation}
If $f$ is a Littlewood polynomial and we want these expressions to be squares (or minus a square -- see Lemma \ref{lem:discrit}), we can count the possible choices of coefficients of $f_{\mr{e}}$ and $f_{\mr{o}}$ giving rise to (possibly degenerate) Pythagorean triples. This is key in the proof of Proposition \ref{prop:nrofsqdisc}.

\begin{proof}[Proof of Proposition \ref{prop:nrofsqdisc}]
Consider a not-necessarily monic reciprocal Littlewood polynomial 
\begin{equation*}
f = a_{4n} X^{8n} + \cdots + a_1 X^{4n+1} + a_0 X^{4n} + a_1 X^{4n-1} + \cdots + a_{4n-1}X + a_{4n}
\end{equation*}
of degree $8n$; since $f$ has square discriminant if and only if $-f$ has square discriminant, we must divide by $2$ whatever final expression we obtain to establish the count of monic reciprocal Littlewood polynomials with square discriminant. Set
\begin{align}
c := f_{\mr{e}}(1) &= a_0 + 2(a_2 + a_4 + \cdots + a_{4n}), \label{eq:c} \\
b := f_{\mr{o}}(1) &= 2(a_1 + a_3 + \cdots + a_{4n-1}). \label{eq:b}
\end{align}
By Lemma \ref{lem:discrit} and \eqref{eq:countrec}, we need to pick the $a_i$ such that $c^2-b^2$ is a square, say equal to $a^2$. In the $\binom{2n}{n}$ cases that exactly half of the odd-index coefficients $a_1$, $a_3$, \ldots, $a_{4n-1}$ are equal to $1$ and thus $b=0$, we find that any choice of the coefficients $a_0$, $a_2$, \ldots, $a_{4n}$ will make $f$ a Littlewood polynomial with square discriminant. There are in total $2^{2n+1} \binom{2n}{n}$ such polynomials. After dividing by two, this is the first term in $\eqref{eq:nrofsqdisc}$.

Now suppose $b$ is nonzero. Recall that if $a^2+b^2 = c^2$ is a Pythagorean triple and $a$, $b$ and $c$ are positive, then there are \textit{unique} positive integers $k$, $r$ and $s$ such that $c=k(r^2+s^2)$, $b = 2krs$, and $a = k(r^2-s^2)$, and $r > s$ and the numbers $r$ and $s$ are coprime and of opposite parity. Since $c$ is odd by definition, we must add the condition that $k$ be odd. This gives the summation condition in \eqref{eq:nrofsqdisc}. The prefactor of $2$ before the sum arises because we treat each of the \textit{four} triples $(a, \pm b, \pm c)$ separately -- we care if $c^2-b^2$ is a square, so the sign of $a$ doesn't matter; but the polynomials corresponding to the four tuples $(\pm b, \pm c)$ are genuinely different. We conclude that the final expression must be multiplied by $4/2 = 2$.

It remains to show that the second summand in \eqref{eq:nrofsqdisc} is correct. That is, we must count all choices of the $a_i$ that lead to the equalities $c = k(r^2+s^2)$ and $b = 2krs$. Notice that
\begin{equation} \label{eq:r^2s^2}
a_2 + a_4 + \cdots + a_{4n} = \frac{c-a_0}{2} = \frac{k(r^2+s^2)-a_0}{2}.
\end{equation}
Since all $a_i$ lie in $\{\pm 1\}$, the left-hand side in \eqref{eq:r^2s^2} is even. As $r^2+s^2 \equiv 1 \bmod{4}$, we find that $a_0 \equiv k \bmod{4}$. Hence $a_0 = (-1)^{\frac{k-1}{2}}$. 
Therefore a total of $n+(k(r^2+s^2)+(-1)^{\frac{k+1}{2}})/4$ of the even-index coefficients $a_2, a_4, \ldots, a_{4n}$ must be equal to $1$. This yields
\begin{equation*}
\binom{2n}{n+\frac{1}{4}(kr^2+ks^2+(-1)^{\frac{k+1}{2}})}
\end{equation*}
options for the even-index coefficients. Similarly, there are $2n$ choices to be made for the odd-index coefficients $a_1, a_3, \ldots, a_{4n-1}$; since the sum of the latter equals $b/2 = krs$, we find that $n+krs/2$ of the odd-index coefficients must be equal to $1$. So we have in total $\binom{2n}{n+krs/2}$ options for the odd-index coefficients. This gives
\begin{equation*}
\binom{2n}{n+\frac{1}{2}krs} \binom{2n}{n+\frac{1}{4}(kr^2+ks^2 + (-1)^{\frac{k+1}{2}})}
\end{equation*}
combinations in total, which is the summand in \eqref{eq:nrofsqdisc}.
\end{proof}

It is clear that the proof method can in principle be applied to derive a combinatorial expression for the number of square-discriminant (skew-)reciprocal polynomials of given degree with coefficients in any fixed set $\mathscr{N}$. For $\abs{R_{8n-2}}$, $\abs{S_{8n}}$ and $\abs{S_{8n-2}}$, we obtain the following expressions.

\begin{lemma} \label{prop:nrofsqdiscother}
We have
\begin{gather}
\abs{R_{8n-2}} =
2^{2n-1} \binom{2n}{n} + 2 \sum \binom{2n}{n+\frac{1}{2}krs} \binom{2n-1}{n+\frac{1}{4}(kr^2+ks^2 + (-1)^{\frac{k-1}{2}}-2)}, \label{eq:R8n-2} \\
\abs{S_{8n}} =
2^{2n} \binom{2n}{n} + 2 \sum \binom{2n}{n+\frac{1}{2}krs} \binom{2n}{n+\frac{1}{4}(kr^2-ks^2 + (-1)^{\frac{k+1}{2}+s})}, \label{eq:S8n} \\
\abs{S_{8n-2}} =
2^{2n-1} \binom{2n}{n} + 2 \sum \binom{2n}{n+\frac{1}{2}krs} \binom{2n-1}{n+\frac{1}{4}(kr^2-ks^2 + (-1)^{\frac{k-1}{2} + s}-2)}, \label{eq:S8n-2}
\end{gather}
where in each case the sum extends over all tuples $(k,r,s)$ such that $k>0$ is odd and $r>s>0$ are coprime and of opposite parity. \qed
\end{lemma}

\section{Lattice point counting}
\label{sec:aux}

The chief aim of this section is to provide integer lattice point count estimates of the regions over which the sums in \eqref{eq:nrofsqdisc} and \eqref{eq:S8n} extend. Define the corresponding sets $B_n$ and $D_n$ as
\begin{gather}
B_n = \left\{(k,r,s) \in \Z^3 \middle\vert \begin{array}{c} k > 0 \text{ and odd, } r > s > 0 \text{ coprime and} \\ \text{of opposite parity, and } k(r^2+s^2) \leq n \end{array} \right\}, \label{eq:Bn} \\
D_n = \left\{(k,r,s) \in \Z^3 \middle\vert \begin{array}{c} k > 0 \text{ and odd, } r > s > 0 \text{ coprime and of} \\ \text{opposite parity, and } k(r^2-s^2) \leq n \text{ and } 2krs \leq n \end{array} \right\}. \label{eq:Dn}
\end{gather}
Since the inequality $k(r^2+s^2) \geq 2krs$ holds for all positive integers $k$, $r$ and $s$, the set $B_{5n}$ certainly contains all tuples $(k,r,s)$ over which the sum in \eqref{eq:nrofsqdisc} extends. The tuples $(k,r,s)$ over which the sum in \eqref{eq:S8n} extends are contained in $D_{5n}$. 

We prove the following asymptotics for the sizes of $B_n$ and $D_n$. These are essentially reproductions of results by Sierpiński \cite[Chapter 5,~Eq.~(7)]{Sierpinski}, and by Benito and Varona \cite[Corollary~2]{BVPyth}, respectively. These sources also specify the error term -- subsequent improvements in the former case are due to Stronina \cite{Stronina} and Nowak and Recknagel \cite{NR}. The differences between the sets $B_n$ and $D_n$ and their equivalents in \cite{Sierpinski} and \cite{BVPyth} are very minor: here, we impose the additional requirements that $k$ be odd and $r>s>0$ in \eqref{eq:Bn}, and that $k$ be odd and $r>s$ in \eqref{eq:Dn}. The proofs given here allow us to derive more specific results, counting such lattice points in certain circle and hyperbolic sectors; see Lemma \ref{lem:latticecount} and Lemma \ref{lem:srecsummary2}. These results are required in the subsequent sections.

\begin{proposition} \label{prop:Bnsize}
The set $B_n$ is of size asymptotically equal to $\frac{1}{4 \pi} n \log{n}$.
\end{proposition}

\begin{proposition} \label{prop:Dnsize}
The set $D_n$ is of size asymptotically equal to $\frac{2\alpha}{\pi^2} n \log{n}$, where $\alpha = \log{\sqrt{1+\sqrt{2}}}$.
\end{proposition}

Let $(a,b,c)$ be a Pythagorean triple, i.e. $a^2+b^2 = c^2$. Assume $a$, $b$ and $c$ are positive integers and $a$ and $b$ are of opposite parity. Proposition \ref{prop:Bnsize} implies that the number of Pythagorean triples with hypotenuse less than $n$ and of opposite parity (considering the triples $(a,b,c)$ and $(b,a,c)$ to be the same) is asymptotic to $\frac{1}{4\pi} n \log{n}$. 
Similarly, Proposition \ref{prop:Dnsize} implies that the number of such triples with legs less than $n$ and of opposite parity is asymptotic to $\frac{2\alpha}{\pi^2} n \log{n}$.

The following corollary of \cite[Theorem 2]{Nymann} is a variant of M\"{o}bius inversion that will be used in the proofs of both Proposition~\ref{prop:Bnsize} and Proposition~\ref{prop:Dnsize}.

\begin{lemma} \label{lem:mobius}
Let $F$ and $f$ be real-valued functions defined on $\R_{\geq 1}$ and related through
\begin{equation*}
F(n) = \sum_{\substack{1 \leq d \leq n \\ d \text{ odd}}} f(n/d).
\end{equation*}
(More explicitly, the sum extends over all odd integers between $1$ and $n$, not just the odd divisors of $n$.) Denoting by $\mu$ the M\"obius function, we have
\begin{equation*}
f(n) = \sum_{\substack{1 \leq d \leq n \\ d \text{ odd}}} \mu(d) F(n/d).
\end{equation*}
\end{lemma}

\subsection*{Reciprocals}
\label{subsec:auxrec}

In this subsection, we prove Proposition~\ref{prop:Bnsize}. Define the circle sector $C_{\theta}(n)$ for $0 \leq \theta \leq \pi/4$ as
\begin{equation*}
C_{\theta}(n) = \{ (x,y) \in \R_{>0}^2 \mid x^2 + y^2 < n^2, \, y \leq x \tan(\theta) \},
\end{equation*}
i.e. the part of the circle of radius $n$ centered at the origin in $\R^2$ that lies in the upper-right quadrant and is bounded by $y = 0$ and $y = x \tan(\theta)$. Denote by $F_{\theta}(n)$ the number of integral, opposite-parity lattice points in the circle sector $C_{\theta}(n)$. Let $f_{\theta}(n)$ be the number of such that are also coprime. 

\begin{lemma} \label{lem:mobiusrec}
We have
\begin{equation*}
f_{\theta}(n) = \sum_{\substack{1 \leq d \leq n \\ d \text{ odd}}} \mu(d) F_{\theta}(n/d).
\end{equation*}
\end{lemma}
\begin{proof}
If $(x,y)$ is an integral, opposite-parity lattice point in the circle sector $C_{\theta}(n)$ with greatest common divisor $d$, then $(x/d, y/d)$ is a primitive, integral, opposite-parity lattice point in the circle sector $C_{\theta}(n/d)$. The opposite holds as well. Noting that a pair of opposite-parity integers that are both at most $n$ must have odd greatest common divisor at most $n$, we find
\begin{equation*}
F_{\theta}(n) = \sum_{\substack{1 \leq d \leq n \\ d \text{ odd}}} f_{\theta}(n/d).
\end{equation*}
Lemma \ref{lem:mobius} gives the desired result.
\end{proof}

The next lemma shows that $F_{\theta}(n)$ and $f_{\theta}(n)$ are linear in $\theta$.

\begin{lemma} \label{lem:latticecount}
The following asymptotics for $F_{\theta}$ and $f_{\theta}$ hold as $n$ goes to infinity:
\begin{enumerate}
    \item $F_{\theta}(n) \sim \theta n^2/4$.
    \item $f_{\theta}(n) \sim 2\theta n^2/\pi^2$.
\end{enumerate}
\end{lemma}
\begin{proof}
Part (b) follows after combining part (a) with Lemma \ref{lem:mobiusrec} and
\begin{equation*}
\sum_{\substack{d \geq 1 \\ d \text{ odd}}} \frac{\mu(d)}{d^2} = \sum_{d \geq 1} \frac{\mu(d)}{d^2} - \sum_{\substack{d \geq 1 \\ d \text{ even}}} \frac{\mu(d)}{d^2} = \sum_{d \geq 1} \frac{\mu(d)}{d^2} - \sum_{d \geq 1} \frac{\mu(2d)}{4d^2} = \sum_{d \geq 1} \frac{\mu(d)}{d^2} + \frac{1}{4} \sum_{\substack{d \geq 1 \\ d \text{ odd}}} \frac{\mu(d)}{d^2}
\end{equation*}
by multiplicativity of the M\"{o}bius function, so that
\begin{equation*}
\sum_{\substack{d \geq 1 \\ d \text{ odd}}} \frac{\mu(d)}{d^2} = \frac{4}{3} \sum_{d \geq 1} \frac{\mu(d)}{d^2} = \frac{8}{\pi^2},
\end{equation*}
see \cite[Corollary 1.10]{MV}.

For part (a), we start by distributing the lattice points in $C_{\theta}(n)$ over four subsets depending on the parity of each of the coordinates. Denote by $F_{00}$ the number of lattice points in $C_{\theta}(n)$ whose coordinates are both even, by $F_{01}$ the number of those whose $x$-coordinate is even and $y$-coordinate is odd, and similarly for $F_{10}$ and $F_{11}$. For each even number $x_0$, the number of lattice points $(x_0, y)$ with odd $y$ exceeds those with even $y$ by at most one. Since $x_0$ lies between $1$ and $n$, we find that $F_{00} + n/2 \geq F_{01}$. Similarly, we deduce $F_{10} + (n+1)/2 \geq F_{11}$ and $F_{00} + n \sin(\theta)/2 \geq F_{10}$. Therefore the difference between any two of the sets $F_{00}$, $F_{01}$, $F_{10}$, and $F_{11}$ is of order $n$.
On the other hand, the quantity $F_{00}$ equals the number of total lattice points in $C_{\theta}(n/2)$. This is asymptotically equal to the area of $C_{\theta}(n/2)$, which is $\theta n^2/8$, see e.g. \cite[Chapter~1.1]{Kratzel}. As $F_{00}$, $F_{01}$, $F_{10}$, and $F_{11}$ differ by a term of order $n$ at most, they are equal asymptotically. Therefore $F_{\theta}(n) = F_{01} + F_{10} \sim \theta n^2/4$.
\end{proof}

We are now ready to prove Proposition~\ref{prop:Bnsize}.

\begin{proof}[Proof of Proposition~\ref{prop:Bnsize}]
Take $\theta = \pi/4$ in Lemma \ref{lem:latticecount}. Then
\begin{equation} \label{eq:Bnsize}
\abs{B_n} = \sum_{\substack{1 \leq d \leq n \\ d \text{ odd}}} f_{\theta}\left(\sqrt{\frac{n}{d}}\right).
\end{equation}
Evaluating the sum up to $d = n^+ := n/\log \log{n}$, we find
\begin{equation*}
\sum_{\substack{1 \leq d \leq n^+ \\ d \text{ odd}}} f_{\theta}\left(\sqrt{\frac{n}{d}}\right) \sim \frac{2 \theta n}{\pi^2} \sum_{\substack{1 \leq d \leq n^+\\ d \text{ odd}}} \frac{1}{d} \sim \frac{\theta n \log{n}}{\pi^2} = \frac{n \log{n}}{4 \pi}
\end{equation*}
where the last asymptotic equality follows since $\log(n/\log \log{n}) \sim \log{n}$. This gives the claimed asymptotic size of $B_n$. The remaining terms of the sum in \eqref{eq:Bnsize}, where $d > n/\log \log{n}$, are bounded by
\begin{equation*}
\sum_{\substack{n^+ < d \leq n \\ d \text{ odd}}} f_{\theta}\left(\sqrt{\frac{n}{d}}\right) < n f_{\theta}\left(\sqrt{\frac{n}{n/\log \log n}}\right) \sim \frac{n \log \log{n}}{2\pi}
\end{equation*}
and thus do not contribute to the asymptotic size of $B_n$.
\end{proof}

\subsection*{Skew-reciprocals}
\label{subsec:auxsrec}

Here, we prove Proposition~\ref{prop:Dnsize} by first establishing skew-reciprocal versions of Lemma \ref{lem:mobiusrec} and Lemma \ref{lem:latticecount}. The situation here is a bit different from the reciprocal case, because neither of $k(r^2-s^2)$ and $2krs$ in the definition of the set $D_n$, see \eqref{eq:Dn}, dominates the other for every choice of positive integers $k$, $r$ and $s$ with $r>s$. Indeed, the inequality $k(r^2-s^2) > 2krs$ holds if and only if $(\sqrt{2}-1)r > s$. Hence both of the inequalities $(r^2-s^2) \leq n$ and $2krs \leq n$ are required in the definition of $D_n$, as opposed to the single inequality $k(r^2+s^2) \leq n$ appearing in the definition of $B_n$, see \eqref{eq:Bn}.

Set $\alpha = \mr{artanh}(\sqrt{2}-1) = \log{\sqrt{1+\sqrt{2}}}$; this is the inverse hyperbolic tangent of the angle between the $r$-axis and the line from the origin to the intersection point of the hyperbolas $r^2-s^2 = n$ and $2rs = n$. Define the hyperbolic sectors $H_{\theta}(n)$ and $H_{\theta}^{*}(n)$ for $0 < \theta \leq \alpha$ as
\begin{align}
H_{\theta}(n) = \{ (x_1,y_1) \in \R^2_{>0} ~|~ x_1^2-y_1^2 < n^2, \, y_1 \leq \tanh(\theta) x_1 \}, \label{eq:Htheta} \\
H_{\theta}^*(n) = \{ (x_2, y_2) \in \R^2_{>0} ~|~ 2x_2y_2 < n^2, \, y_2 < x_2 \leq e^{2\theta} y_2 \}. \label{eq:Hstartheta}
\end{align}
Note that $e^{2\theta} = (1+\tanh(\theta))/(1-\tanh(\theta))$, and that both $\tanh(\alpha)$ and $e^{2\alpha}$ are equal to $\sqrt{2}-1$. 

As in the reciprocal case, the reason to consider these sectors is that their areas scale linearly in $\theta$.
\begin{lemma} \label{lem:areahypsectors}
The areas of $H_{\theta}(n)$ and of $H^*_{\theta}(n)$ are each equal to $\theta n^2/2$.
\end{lemma}
\begin{proof}
The linear transformation sending $x_2 \mapsto (x_1+y_1)/\sqrt{2}$ and $y_2 \mapsto (x_1-y_1)/\sqrt{2}$ maps $H^*_{\theta}(n)$ to $H_{\theta}(n)$ and has determinant $1$. Thus $H^*_{\theta}(n)$ and $H_{\theta}(n)$ have equal area. The area of $H_{\theta}(n)$ is $n^2$ times as large as that of the region bounded by the hyperbola $x_1^2-y_1^2 = 1$, the axis $y_1 = 0$, and the ray through the origin and the point $(\cosh(\theta), \sinh(\theta))$. But that is simply $\theta/2$.
\end{proof}

We summarise the analogues of Lemma \ref{lem:mobiusrec} and Lemma \ref{lem:latticecount} in the following lemma.

\begin{lemma} \label{lem:srecsummary2}
Denote by $G_{\theta}(n)$ (resp. $G^*_{\theta}(n)$) the number of integral, opposite-parity lattice points in $H_{\theta}(n)$ (resp. $H^*_{\theta}(n)$), and by $g_{\theta}(n)$ (resp. $g^*_{\theta}(n)$) the number of such that are also coprime. Then the following hold:
\begin{enumerate}
    \item $g(n) = \sum \mu(d) G(n/d)$ where the sum extends over all odd $1 \leq d \leq n$, and similarly for $g^*_{\theta}$.
    \item $G_{\theta}(n) \sim G^*_{\theta}(n) \sim \theta n^2/4$.
    \item $g_{\theta}(n) \sim g^*_{\theta}(n) \sim 2\theta n^2/\pi^2$.
\end{enumerate}
\end{lemma}
\begin{proof}[Proof sketch]
All proofs are analogous to those of the mentioned lemmas, where $H_{\theta}(n)$ (respectively $H^*_{\theta}(n)$) plays the role of $C_{\theta}(n)$. That the asymptotic expressions for $f_{\theta}$, $g_{\theta}$, and $g^*_{\theta}$ are all equal comes from the fact that the circle sector $C_{\theta}$ and the hyperbolic sectors $H_{\theta}$ and $H^*_{\theta}$ all have equal area, see Lemma \ref{lem:areahypsectors}.
\end{proof}

We are now in the position to prove Proposition~\ref{prop:Dnsize}.

\begin{proof}[Proof of Proposition~\ref{prop:Dnsize}]
Note that
\begin{equation} \label{eq:Dnsize}
\abs{D_n} = \sum_{\substack{1 \leq d \leq n \\ d \text{ odd}}} g_{\alpha}\left( \sqrt{\frac{n}{d}} \right) + g^*_{\alpha}\left( \sqrt{\frac{n}{d}} \right).
\end{equation}
Writing $n^{+} = n/\log \log{n}$ and reasoning as in the proof of Proposition \ref{prop:Bnsize} that the terms in the sum with $d > n^+$ do not contribute, we find with help of Lemma \ref{lem:srecsummary2}(c) that
\begin{equation*}
\abs{D_n} \sim \sum_{\substack{1 \leq d \leq n^{+} \\ d \text{ odd}}} g_{\alpha}\left( \sqrt{\frac{n}{d}} \right) + g^*_{\alpha}\left( \sqrt{\frac{n}{d}} \right) \sim \frac{4 \alpha n}{\pi^2} \sum_{\substack{ 1 \leq d \leq n^+ \\ d \text{ odd}}} \frac{1}{d} \sim \frac{2 \alpha n \log{n}}{\pi^2}
\end{equation*}
as claimed.
\end{proof}

\section{The reciprocals}
\label{sec:reciprocals}

In this section, we build up towards the proof of the part of Theorem \ref{thm:mainthm} that concerns reciprocals. For the proof, we break up the sum in \eqref{eq:nrofsqdisc} into several pieces. Fix a (large) integer $N$ and set $\epsilon = N^{-1}$ and $m = 5\sqrt{n \log{n}}$ (the number $5$ is a convenient choice, but could be replaced by any real number greater than $2\sqrt{2}$). Write
\begin{align*}
\Sigma_1 &= \sum_{(k,r,s) \in B_{N\sqrt{n}}} \binom{2n}{n+\frac{1}{2}krs} \binom{2n}{n+\frac{1}{4}(kr^2+ks^2 + (-1)^{\frac{k+1}{2}})}, \\
\Sigma_2 &= \sum_{(k,r,s) \in B_m \setminus B_{N\sqrt{n}}} \binom{2n}{n+\frac{1}{2}krs} \binom{2n}{n+\frac{1}{4}(kr^2+ks^2 + (-1)^{\frac{k+1}{2}})},
\end{align*}
and define $\Sigma_3$ through $2\Sigma_3 = \abs{R_{8n}} - 2^{2n}\binom{2n}{n} - 2\Sigma_1 - 2\Sigma_2$. Figure \ref{fig:1} shows how the domain consisting of lattice points over which the sum in \eqref{eq:nrofsqdisc} extends is divided into parts associated with the sums $\Sigma_1$, $\Sigma_2$ and $\Sigma_3$. The following subsections go into the asymptotics of each of these terms, showing that $\Sigma_1$ is the dominant term. To obtain an exact expression for the main term in the asymptotics of $\Sigma_1$, precise control over both binomial coefficients in its summand is needed. In contrast, to show that $\Sigma_2$ is negligible in comparison, we only need to control one binomial coefficient precisely, and for $\Sigma_3$ it suffices to estimate both binomial coefficients appearing in the summand by the maximal value they can obtain. Proposition \ref{prop:almostcentral} and Proposition \ref{prop:Bnsize} are key in this.

\begin{figure}[ht]
    \centering
    \begin{tikzpicture}
    \draw [<->,thick] (0,4.5) node (yaxis) [above] {$s$}
        |- (6.8,0) node (xaxis) [right] {$r$};
\draw[very thick] (0,0) -- (45:6) arc(45:0:6) node[midway,left]{$\Sigma_3$} node[below]{$\sqrt{\frac{5n}{k}}$} -- cycle;
\draw[very thick,fill=black!25] (0,0) -- (45:4.7) arc(45:0:4.7) node[midway,left]{$\Sigma_2$} node[below]{$\sqrt{\strut\frac{m}{k}}$} -- cycle;
\draw[very thick,fill=black!50] (0,0) -- (45:3) arc(45:0:3) node[midway,left]{$\Sigma_1$} node[below]{$\sqrt{\frac{N\sqrt{n}}{k}}$} -- cycle; 
\draw (45:6) -- (45:6.5) node[above]{$r=s$};
\draw[latex-latex]  (45:1.2) arc(45:0:1.2) node[midway,right]{$\frac{\pi}{4}$};
\end{tikzpicture}
    \caption{Slice of the domain containing $B_{5n}$ at a fixed $k$, showing the subdomains related to the sums $\Sigma_i$ with $i = 1,2,3$. The full domain (with $k$ varying) is part of the interior of an elliptic paraboloid.}
    \label{fig:1}
\end{figure}

We often use elementary estimates of sums by integrals without reference; proofs for any such estimate may be found in \cite[Theorems 4.1 and 4.2]{Spencer}. In addition, here we also record the following asymptotic of binomial coefficients that are close to central (although we will also use it for the skew-reciprocals). We refer to \cite[\S 5.4]{Spencer} for a proof.

\begin{proposition} \label{prop:almostcentral}
We have $\binom{2n}{n} \sim 4^n/\sqrt{\pi n}$. Furthermore, if $k$ is of order $o(n^{2/3})$, then
\begin{equation*}
\binom{2n}{n+k} \sim \binom{2n}{n} e^{-\frac{k^2}{n}}.
\end{equation*}
\end{proposition}

\subsection*{The sums \texorpdfstring{$\Sigma_2$}{Sigma2} and \texorpdfstring{$\Sigma_3$}{Sigma3}}
\label{subsec:S2S3}

In this subsection, we show that the term $2^{2n} \binom{2n}{n}$ and the sums $\Sigma_2$ and $\Sigma_3$ each have negligible contribution in comparison to $16^n \log{n}/\sqrt{n}$ when $\epsilon$ tends to zero. First, recall that we have already seen in the introduction that $2^{2n} \binom{2n}{n} \asymp 16^n/\sqrt{n}$. The sum $\Sigma_3$ satisfies
\begin{equation} \label{eq:S3}
\Sigma_3 \leq \abs{B_{5n}} \binom{2n}{n} \binom{2n}{n+\sqrt{n \log n}} \asymp \abs{B_{5n}} \binom{2n}{n}^2 e^{-\frac{n \log{n}}{n}} \asymp \frac{16^n \log{n}}{n}
\end{equation}
by Proposition \ref{prop:Bnsize} and Proposition \ref{prop:almostcentral}. In conclusion, both $2^{2n}\binom{2n}{n}$ and $\Sigma_3$ are of order $o(16^n \log{n}/\sqrt{n})$. 

\begin{lemma} \label{lem:S2}
The sum $\Sigma_2$ satisfies
\begin{equation*}
\lim_{\epsilon \to 0} \lim_{n \to \infty} \frac{\Sigma_2}{16^n \log{n}/\sqrt{n}} = 0.
\end{equation*}
\end{lemma}
\begin{proof}
Defining $C = C(k,r,s) = 4k^2r^2s^2 + (k(r^2+s^2) + (-1)^{\frac{k+1}{2}})^2$, Proposition \ref{prop:almostcentral} implies that
\begin{equation*}
\Sigma_2 \sim \binom{2n}{n}^{2} \sum_{\substack{(k,r,s) \in B_m \setminus B_{N\sqrt{n}}}} e^{-\frac{C}{16n}}.
\end{equation*}
Since $C > k^2(r^2+s^2)^2$ for all positive integers $k$, $r$, and $s$, the sum on the right-hand side is bounded from above by
\begin{equation*}
\sum_{\substack{1 \leq r \leq \sqrt{m} \\ 2 \leq s \leq \sqrt{m} \\ k > \frac{N\sqrt{n}}{r^2+s^2}}} e^{-\frac{1}{16n}k^2(r^2+s^2)^2} \leq \int_1^{\sqrt{m}} \int_0^{\sqrt{m}} \left( 1 + \int_{\frac{N\sqrt{n}}{r^2+s^2}}^{\infty} e^{-\frac{1}{16n}k^2(r^2+s^2)^2} \dif k \right) \dif r \dif s
\end{equation*}
by applying elementary estimates for sums by integrals. Pulling out the $1$ from the middle integral and evaluating the innermost integral yields
\begin{equation*}
\sum_{\substack{(k,r,s) \in B_m \setminus B_{N\sqrt{n}}}} e^{-\frac{C}{16n}} < m + 2\sqrt{\pi n} (1 - \erf(N/4)) \int_1^{\sqrt{m}} \int_0^{\sqrt{m}} \frac{1}{r^2+s^2}  \dif r \dif s,
\end{equation*}
where $\erf(x) = 2\pi^{-1/2} \int_0^x e^{-t^2} \dif t$ is the error function.
Switching to polar coordinates with $R^2 = r^2 + s^2$, the remaining double integral is bounded by
\begin{equation*}
\int_1^{\sqrt{m}} \int_0^{\sqrt{m}} \frac{1}{r^2+s^2}  \dif r \dif s < \frac{\pi}{2} \int_{1}^{\sqrt{2m}} \frac{1}{R} \dif R = \frac{\pi}{4} \log{2m} \asymp \log{n}.
\end{equation*}
Thus, as $\binom{2n}{n}^2 \asymp 16^n/n$, the sum $\Sigma_2$ is asymptotically at most
\begin{equation*}
(1 - \erf(N/4)) \frac{16^n \log{n}}{\sqrt{n}}
\end{equation*}
up to a multiplicative constant independent of $n$ and $N$. As $\erf(x)$ goes to $1$ as $x$ tends to $\infty$, this yields the claimed limit.
\end{proof}

\subsection*{The sum \texorpdfstring{$\Sigma_1$}{Sigma1}}
\label{subsec:S1}

\begin{figure}[ht]
    \centering
    \begin{tikzpicture}
    \draw [<->,thick] (0,4.5) node (yaxis) [above] {$s$}
        |- (6.8,0) node (xaxis) [right] {$r$};
\draw[very thick,fill=black!50] (0,0) -- (45:6) arc(45:0:6) node[below]{$\sqrt{\frac{N\sqrt{n}}{k}}$} -- cycle; 
\draw (45:1) arc(45:0:1);
\draw (45:2) arc(45:0:2) node[below]{$\ldots$};
\draw (45:3) arc(45:0:3) node[below]{$t_{j-1}$};
\draw (45:4) arc(45:0:4) node[below]{$t_j$};
\draw (45:5) arc(45:0:5) node[below]{$\ldots$};
\draw[fill=black!20] (27:3) coordinate (beta) arc (27:18:3) coordinate (alpha) -- (18:4) arc (18:27:4) -- cycle;
\draw (3.25,1.33) node{$T_{ij}$};
\draw (4.2,3.5) node{$\Sigma_1$};
\draw (1.7, 0.71) node{$\theta$};
\draw (0,0) -- (9:6) arc(9:0:6);
\draw (0,0) -- (18:6) arc(18:0:6);
\draw (0,0) -- (27:6) arc(27:0:6);
\draw (0,0) -- (36:6) arc(36:0:6);
\draw (45:6) -- (45:6.5) node[above]{$r=s$};
\draw (9:6) -- (9:6.2) node[right]{$\ldots$};
\draw (18:6) -- (18:6.2) node[right]{$s/r = \tan(\theta_{i-1})$};
\draw (27:6) -- (27:6.2) node[right]{$s/r = \tan(\theta_i)$};
\draw (36:6) -- (36:6.2) node[right]{$\ldots$};
\draw[latex-latex]  (27:1.6) arc(27:18:1.6);
\end{tikzpicture}
    \caption{Slice of the domain over which the sum $\Sigma_1$ extends at a fixed $k$, illustrating a two-dimensional section of the radial grid defined by the inequalities \eqref{eq:ineq1} and \eqref{eq:ineq2}. Here, we have $t_{\ell} := \sqrt{\frac{\ell \epsilon \sqrt{n}}{k}}$ and $\theta := \epsilon \frac{\pi}{4}$. The light grey-shaded grid cell labelled $T_{ij}$ contains, by abuse of notation, a two-dimensional slice of the lattice point set $T_{i,j} \subset \Z^3$ of the same name, see \eqref{eq:Tij}.}
    \label{fig:2}
\end{figure}

To obtain a precise estimate of $\Sigma_1$, we need to control both binomial coefficients in the summand of \eqref{eq:nrofsqdisc} simultaneously. This is achieved by dividing the domain over which the sum extends in boxes as follows. Let $1 \leq i \leq N$ and $1 \leq j \leq N^2$. Write $\theta_i = i\epsilon\pi/4$ and consider the inequalities
\begin{gather}
(j-1)\epsilon\sqrt{n} < k(r^2+s^2) \leq j\epsilon \sqrt{n}, \label{eq:ineq1} \\
\tan(\theta_{i-1}) < s/r \leq \tan(\theta_i);\label{eq:ineq2}
\end{gather}
this is a region enclosed between two circles and two lines. For fixed positive $k$, the inequalities \eqref{eq:ineq1} and \eqref{eq:ineq2} partition half of the right-upper quadrant of the disk $r^2 + s^2 \leq N\sqrt{n}/k$ in a radial grid, see Figure~\ref{fig:2}. 
Define the set $T_{ij}$ as
\begin{equation} \label{eq:Tij}
T_{ij} = \left\{(k,r,s) \in \Z^3 \, \middle\vert \begin{array}{c} k > 0 \text{ and odd, } r > s > 0 \text{ coprime and of opposite} \\ \text{parity, and } (k,r,s) \text{ satisfies } \eqref{eq:ineq1} \text{ and } \eqref{eq:ineq2} \end{array} \right\}
\end{equation}
-- see Figure~\ref{fig:2} for an illustration. The following lemma demonstrates that the lattice point sets $T_{ij}$ are asymptotically equal in size.

\begin{lemma} \label{lem:Tij}
As $n$ tends to infinity, we have
\begin{equation*}
\abs{T_{ij}} \sim \frac{\epsilon^2}{8\pi} \sqrt{n} \log{n}.
\end{equation*}
In particular, the size of $T_{ij}$ does not depend on $i$ and $j$.
\end{lemma}

\begin{proof}
Write $a = \sqrt{j\epsilon \sqrt{n}/k}$ and $b = \sqrt{(j-1)\epsilon \sqrt{n}/k}$. For fixed positive $k$, the number of integral, coprime, opposite-parity lattice points $(r, s)$ in the box bounded by the inequalities \eqref{eq:ineq1} and \eqref{eq:ineq2} equals
\begin{equation}
\label{eq:qij}
q(i,j) := \left(f_{\theta_i}(a) -f_{\theta_{i-1}}(a)\right) - \left(f_{\theta_i}(b) - f_{\theta_{i-1}}(b)\right).
\end{equation}
For $k \leq n^+ := \sqrt{n}/\log{\log{n}}$, the quantity $q$ satisfies the asymptotic equality
\begin{equation*}
q(i,j) \sim \frac{2}{\pi^2} (\theta_i - \theta_{i-1})(a^2-b^2) = \frac{\epsilon^2}{2 \pi} \frac{\sqrt{n}}{k}
\end{equation*}
by Lemma \ref{lem:latticecount} (notice that this is just the area of the box multiplied by $4/
\pi^2$). When $k > n^+$, the bound $q(i,j) \ll \log{\log{n}}$ holds as each of the four terms on the right-hand side of \eqref{eq:qij} are at most of this order.
Therefore
\begin{equation*}
\sum_{\substack{ 1 \leq k \leq n^+ \\ k \text{ odd}}} q(i,j) \sim \frac{\epsilon^2}{2 \pi} \sqrt{n} \sum_{\substack{ 1 \leq k \leq n^+ \\ k \text{ odd}}} \frac{1}{k} \sim \frac{\epsilon^2}{8\pi} \sqrt{n} \log{n}
\quad \text{and} \quad
\sum_{\substack{ n^+ < k \leq \sqrt{n} \\ k \text{ odd}}} q(i,j) \ll \sqrt{n} \log{\log{n}},
\end{equation*}
which implies
\begin{equation*}
\abs{T_{ij}} = \sum_{\substack{ 1 \leq k \leq \sqrt{n} \\ k \text{ odd}}} q(i,j) \sim \frac{\epsilon^2}{8\pi} \sqrt{n} \log{n}
\end{equation*}
as claimed.
\end{proof}

Multiplying \eqref{eq:ineq2} through by $r$ and using both of the resulting inequalities,  some rewriting of \eqref{eq:ineq1} leads to
\begin{equation*}
\frac{1}{4} m(i, \epsilon) (j-1)\epsilon \sqrt{n} < \frac{1}{2}krs < \frac{1}{4} M(i, 
\epsilon) j\epsilon \sqrt{n},
\end{equation*}
where
\begin{equation*}
m(i, \epsilon) = 2\tan(\theta_{i-1})\cos^2(\theta_i) \qquad \text{and} \qquad M(i, \epsilon) = 2\tan(\theta_i)\cos^2(\theta_{i-1}).
\end{equation*}
Note that $m(i,\epsilon)$ is increasing on the interval $[-1, N]$ and $M(i, \epsilon)$ is increasing on $[-1, N+1]$.

\begin{lemma} \label{lem:S1}
The sum $\Sigma_1$ satisfies
\begin{equation} \label{eq:sinint}
\lim_{\epsilon \to 0} \lim_{n \to \infty} \frac{\Sigma_1}{16^n \log{n}/\sqrt{n}} = \frac{1}{4\pi^{3/2}} \int_0^1 \frac{1}{\sqrt{1 + \sin^2(\pi t/2) }} \dif t.
\end{equation}
\end{lemma}
\begin{proof}
We give an upper and a lower bound that converge to the same value as $\epsilon$ tends to $0$. For the upper bound, note that
\begin{align}
\Sigma_1 &= \sum_{\substack{1 \leq i \leq N \\ 1 \leq j \leq N^2}} \sum_{(k,r,s) \in T_{ij}} \binom{2n}{n+\frac{1}{2} krs} \binom{2n}{n + \frac{1}{4} (kr^2+ks^2 + (-1)^{(k+1)/2})} \label{eq:S1sumTij} \\ 
&\leq \binom{2n}{n}^2 \sum_{1 \leq i \leq N} \abs{T_{i1}} + \sum_{\substack{1 \leq i \leq N \\ 2 \leq j \leq N^2}} \abs{T_{ij}} \binom{2n}{n+ \frac{1}{4} m(i, \epsilon) (j-1) \epsilon \sqrt{n}} \binom{2n}{n + \frac{1}{4} (j-1) \epsilon \sqrt{n} - 1}. \nonumber
\end{align}
The first sum in the last line, where $j=1$ is fixed, has negligible contribution as $\epsilon \to 0$. In addition, the asymptotics of the last binomial coefficient is not altered by changing $n + \frac{1}{4} (j-1) \epsilon \sqrt{n} - 1$ to $n + \frac{1}{4} (j-1) \epsilon \sqrt{n}$. Combined with Proposition \ref{prop:almostcentral} and Lemma \ref{lem:Tij}, the sum $\Sigma_1$ is therefore asymptotically no larger than
\begin{equation} \label{eq:bigsum}
\frac{\epsilon^2}{8\pi^2} \frac{16^n \log{n}}{\sqrt{n}} \sum_{i=1}^{N} \sum_{j=2}^{N^2} e^{-\frac{1}{16} (j-1)^2 \epsilon^2 (1+m(i,\epsilon)^2)}.
\end{equation}
The inner sum in \eqref{eq:bigsum} is smaller than
\begin{equation*}
\int_1^{\infty} e^{-\frac{1}{16} (j-1)^2 \epsilon^2 (1+m(i,\epsilon)^2)} \dif j = \frac{2 \sqrt{\pi}}{\epsilon \sqrt{1+m(i,\epsilon)^2}}.
\end{equation*}
Plugging this into \eqref{eq:bigsum}, moving out all constants from the sum but keeping all $\epsilon$'s in it shows that it remains to evaluate
\begin{equation*}
\sum_{i=1}^{N} \frac{\epsilon}{\sqrt{1+m(i,\epsilon)^2}}.
\end{equation*}
Again, we employ an integral estimate (using that $m(i,\epsilon)$ is increasing on the interval $[0, N]$) to bound the last sum from above by
\begin{equation*}
\int_0^{N} \frac{\epsilon}{\sqrt{1+m(i,\epsilon)^2}} \dif i =
\int_0^1 \frac{1}{\sqrt{ 1 + 4 \tan^2\left(\frac{(x-\epsilon) \pi}{4}\right)\cos^4\left(\frac{x\pi}{4}\right)}} \dif x
\end{equation*}
after the substitution $x = i\epsilon$. As $\epsilon$ tends to $0$, this becomes the integral shown in \eqref{eq:sinint}.

Now we prove that the asymptotic lower bound is the same. Starting from \eqref{eq:S1sumTij}, notice that this can be bounded from below by
\begin{equation*}
\sum_{\substack{1 \leq i \leq N \\ 1 \leq j \leq N^2}} \abs{T_{ij}} \binom{2n}{n+ \frac{1}{4} M(i, \epsilon) j \epsilon \sqrt{n}} \binom{2n}{n + \frac{1}{4} j \epsilon \sqrt{n} + 1}.
\end{equation*}
Again, dropping the $+1$ in the last binomial coefficient, this sum is asymptotically at least
\begin{equation*}
\frac{\epsilon^2}{8\pi^2} \frac{16^n \log{n}}{\sqrt{n}} \sum_{i=1}^{N} \sum_{j=1}^{N^2} e^{-\frac{1}{16} j^2 \epsilon^2 (1+M(i,\epsilon)^2)}.
\end{equation*}
The inner sum is at least
\begin{equation*}
\int_1^{N^2} e^{-\frac{1}{16} j^2 \epsilon^2 (1+M(i,\epsilon)^2)} \dif j = \frac{2 \sqrt{\pi}}{\epsilon \sqrt{1 + M(i, \epsilon)^2}}  \left( \erf\left( \frac{\sqrt{1+M(i,\epsilon)^2}}{4\epsilon} \right) - \erf\left( \frac{\epsilon\sqrt{(1+M(i,\epsilon)^2)}}{4} \right) \right).
\end{equation*}
Since the error function is monotonously increasing, and $M(i,\epsilon)$ is monotonously increasing on $[1, N]$ as well, the term involving the error functions is at least
\begin{equation*}
\erf\left( \frac{1}{4\epsilon} \right) - \erf\left( \frac{\epsilon\sqrt{3}}{4} \right)
\end{equation*}
which tends to $1$ as $\epsilon$ tends to $0$. We are left with the sum
\begin{equation*}
\sum_{i=1}^{N} \frac{\epsilon}{\sqrt{1 + M(i, \epsilon)^2}}
\end{equation*}
which is bounded from below by
\begin{equation*}
\int_1^{N} \frac{\epsilon}{\sqrt{1 + M(i, \epsilon)^2}} \dif i =
\int_{\epsilon}^1 \frac{1}{\sqrt{ 1 + 4 \tan^2(\frac{x \pi}{4})\cos^4(\frac{(x-\epsilon)\pi}{4})}} \dif x
\end{equation*}
where $x = i \epsilon$. In the limit $\epsilon \to 0$ this becomes the integral on the right-hand side in \eqref{eq:sinint}.
\end{proof}

\section{The skew-reciprocals}
\label{sec:skewreciprocals}

As in the reciprocal case, fix some (large) integer $N$, define $\epsilon = N^{-1}$ and $m = 5\sqrt{n \log{n}}$, and write
\begin{equation*}
\abs{S_{8n}} = 2^{2n} \binom{2n}{n} + 2\Sigma'_1 + 2\Sigma'_2 + 2\Sigma'_3
\end{equation*}
where
\begin{align*}
\Sigma'_1 &= \sum_{(k,r,s) \in D_{N\sqrt{n}}} \binom{2n}{n+\frac{1}{2}krs} \binom{2n}{n+\frac{1}{4}(kr^2-ks^2 + (-1)^{s+\frac{k+1}{2}})}, \\
\Sigma'_2 &= \sum_{(k,r,s) \in D_m \setminus D_{N\sqrt{n}}} \binom{2n}{n+\frac{1}{2}krs} \binom{2n}{n+\frac{1}{4}(kr^2-ks^2 + (-1)^{s+\frac{k+1}{2}})}.
\end{align*}
With methods very similar to the ones employed in the reciprocal case, in the double limit as first $n$ and then $N$ tends to infinity, each of $2^{2n} \binom{2n}{n}$ and the sums $\Sigma'_2$ and $\Sigma'_3$ are negligible compared to $16^n \log{n}/\sqrt{n}$. Here, we focus on the evaluation of $\Sigma'_1$.

Let $1 \leq i \leq N$ and $1 \leq j \leq N^2$. Recall that $\alpha$ is the constant $\log{\sqrt{1+\sqrt{2}}}$. Write $\theta_i = i \epsilon \alpha$ and consider the inequalities
\begin{gather}
(j-1)\epsilon \sqrt{n} < k(r^2-s^2) \leq j\epsilon \sqrt{n}, \label{eq:ineq3} \\
\tanh(\theta_{i-1}) < s/r \leq \tanh(\theta_i); \label{eq:ineq4}
\end{gather}
this is a region enclosed between two hyperbolas and two lines. The quantity $\tanh(\theta)$ varies between $0$ and $\tanh(\alpha) = \sqrt{2}-1$ as $\theta$ varies between $0$ and $\alpha$. Therefore, fixing $k$, the regions described by the inequalities \eqref{eq:ineq3} and \eqref{eq:ineq4} partition $H_{\alpha}(\sqrt{j\epsilon\sqrt{n}/k})$ (with $H_{\alpha}$ as in \eqref{eq:Htheta}). Similarly, the inequalities
\begin{gather}
(j-1)\epsilon \sqrt{n} < 2krs \leq j\epsilon \sqrt{n}, \label{eq:ineq5} \\
e^{-2\theta_{i}} < s/r \leq e^{-2\theta_{i-1}} \label{eq:ineq6}
\end{gather}
partition $H^*_{\alpha}(\sqrt{j\epsilon\sqrt{n}/k})$. 
Define the set
\begin{equation*}
T'_{ij} = \left\{(k,r,s) \in \Z^3 \, \middle\vert \begin{array}{c} k > 0 \text{ and odd, } r > s > 0 \text{ coprime and of opposite} \\ \text{parity, and } (k,r,s) \text{ satisfies } \eqref{eq:ineq3} \text{ and } \eqref{eq:ineq4} \end{array} \right\},
\end{equation*}
and let $T'^{*}_{ij}$ be the similar set of tuples that satisfy \eqref{eq:ineq5} and \eqref{eq:ineq6} instead.

\begin{lemma} \label{lem:Tij'}
As $n$ tends to infinity, we have
\begin{equation*}
\abs*{T'_{ij}} \sim \abs*{T'^{*}_{ij}} \sim \frac{\alpha\epsilon^2}{2\pi^2} \sqrt{n} \log{n}.
\end{equation*}
In particular, the size of $T'_{ij}$ does not depend on $i$ and $j$.
\end{lemma}
\begin{proof}
We argue as in Lemma \ref{lem:Tij}. Write $a = \sqrt{j\epsilon \sqrt{n}/k}$ and $b = \sqrt{(j-1)\epsilon \sqrt{n}/k}$. For fixed positive $k$, the number of integral, coprime, opposite-parity lattice points in the box bounded by the inequalities \eqref{eq:ineq3} and \eqref{eq:ineq4} equals
\begin{equation} \label{eq:qij2}
q(i,j) := \left(g_{\theta_i}(a) -g_{\theta_{i-1}}(a)\right) - \left(g_{\theta_i}(b) - g_{\theta_{i-1}}(b)\right).
\end{equation}
For $k \leq n^+ := \sqrt{n}/\log{\log{n}}$, we deduce the asymptotic equality
\begin{equation*}
q(i,j) \sim \frac{2}{\pi^2} (\theta_i - \theta_{i-1})(a^2-b^2) = \frac{2\alpha\epsilon^2 }{\pi^2} \frac{\sqrt{n}}{k}
\end{equation*}
by Lemma \ref{lem:srecsummary2}. When $k > n^+$, the bound $q(i,j) \ll \log{\log{n}}$ holds as each of the four terms on the right-hand side in \eqref{eq:qij2} are at most of this order. By an argument entirely similar to the one in Lemma \ref{lem:Tij}, we find
\begin{equation*}
\abs*{T'_{ij}} = \sum_{\substack{ 1 \leq k \leq \sqrt{n} \\ k \text{ odd}}} q(i,j) \sim \sum_{\substack{ 1 \leq k \leq n^+ \\ k \text{ odd}}} q(i,j) \sim \frac{2\alpha\epsilon^2}{\pi^2} \sqrt{n} \sum_{\substack{ 1 \leq k \leq n^+ \\ k \text{ odd}}} \frac{1}{k} \sim \frac{\alpha\epsilon^2}{2\pi^2} \sqrt{n} \log{n},
\end{equation*}
as claimed. The same argument gives the result for $T'^{*}_{ij}$.
\end{proof}

Write
\begin{equation*}
m'(i, \epsilon) = \sinh(2\theta_{i-1}) \qquad \text{and} \qquad M'(i, \epsilon) = \sinh(2\theta_{i}).
\end{equation*}
Manipulating the inequalities \eqref{eq:ineq3} and \eqref{eq:ineq4} leads to
\begin{equation*}
\frac{1}{4} m'(i, \epsilon) (j-1)\epsilon \sqrt{n} < \frac{1}{2}krs \leq \frac{1}{4} M'(i, \epsilon) j\epsilon \sqrt{n},
\end{equation*}
whereas the inequalities \eqref{eq:ineq5} and \eqref{eq:ineq6} yield
\begin{equation*}
\frac{1}{4}m'(i, \epsilon) (j-1)\epsilon \sqrt{n} < \frac{1}{4}k(r^2-s^2) \leq \frac{1}{4}M'(i, \epsilon) j\epsilon \sqrt{n}
\end{equation*}
for the same functions $m'$ and $M'$. Write
\begin{equation} \label{eq:S1sumTijS}
\Sigma = \sum_{\substack{1 \leq i \leq N \\ 1 \leq j \leq N^2}} \sum_{(k,r,s) \in T'_{ij}} \binom{2n}{n+\frac{1}{2}krs} \binom{2n}{n+\frac{1}{4}(kr^2-ks^2 + (-1)^{s+\frac{k+1}{2}})}
\end{equation}
and
\begin{equation*}
\Sigma^* = \sum_{\substack{1 \leq i \leq N \\ 1 \leq j \leq N^2}} \sum_{(k,r,s) \in T'^*_{ij}} \binom{2n}{n+\frac{1}{2}krs} \binom{2n}{n+\frac{1}{4}(kr^2-ks^2 + (-1)^{s+\frac{k+1}{2}})},
\end{equation*}
so that $\Sigma_1' = \Sigma + \Sigma^*$.
\begin{lemma}
Each of the sums $\Sigma$ and $\Sigma^*$ can be asymptotically bounded from above by
\begin{equation} \label{eq:bigsumSabove}
\frac{\alpha\epsilon^2}{2\pi^3} \frac{16^n \log{n}}{\sqrt{n}} \sum_{i=1}^{N} \sum_{j=1}^{N^2} e^{-\frac{1}{16} (j-1)^2 \epsilon^2 (1+m'(i,\epsilon)^2)}
\end{equation}
and from below by
\begin{equation} \label{eq:bigsumSbelow}
\frac{\alpha\epsilon^2}{2\pi^3} \frac{16^n \log{n}}{\sqrt{n}} \sum_{i=1}^{N} \sum_{j=1}^{N^2} e^{-\frac{1}{16} j^2 \epsilon^2 (1+M'(i,\epsilon)^2)}.
\end{equation}
In particular, $\Sigma_1'$ is asymptotically equal to $2 \Sigma$.
\end{lemma}
\begin{proof}
We give an upper and a lower bound that converge to the same value as $\epsilon$ tends to $0$. For the upper bound, note that
\begin{equation*}
\Sigma \leq \binom{2n}{n}^2 \sum_{1 \leq i \leq N} \abs*{T'_{i1}} + \sum_{\substack{1 \leq i \leq N \\ 2 \leq j \leq N^2}} \abs*{T'_{ij}} \binom{2n}{n+ \frac{1}{4} m'(i, \epsilon) (j-1) \epsilon \sqrt{n}} \binom{2n}{n + \frac{1}{4} (j-1) \epsilon \sqrt{n} - 1}.
\end{equation*}
The $-1$ appearing in the last binomial coefficient can simply be ignored, because it doesn't affect the asymptotics in $n$ of that binomial coefficient. In addition, we see that the first term in the last line will have negligible contribution as $\epsilon \to 0$. The asymptotics for almost central binomial coefficients given in Proposition \ref{prop:almostcentral} and for $\abs*{T_{ij}'}$ of Lemma \ref{lem:Tij'} show the last sum is asymptotically no larger than the sum in \eqref{eq:bigsumSabove}.

For the lower bound, we observe
\begin{equation*}
\Sigma \geq \sum_{\substack{1 \leq i \leq N \\ 1 \leq j \leq N^2}} \abs*{T'_{ij}} \binom{2n}{n+ \frac{1}{4} M'(i, \epsilon) j \epsilon \sqrt{n}} \binom{2n}{n + \frac{1}{4} j \epsilon \sqrt{n} + 1}
\end{equation*}
starting from \eqref{eq:S1sumTijS}. Again, dropping the $+1$ in the last binomial coefficient, this sum is asymptotically at least the sum in \eqref{eq:bigsumSbelow}.
After replacing $T'_{ij}$ by $T'^*_{ij}$, the same argument holds for $\Sigma^*$.
\end{proof}

We are now in the position to obtain our main result for $\Sigma_1'$.
\begin{lemma}
With $\alpha = \log{\sqrt{1+\sqrt{2}}}$, the sum $\Sigma_1'$ satisfies
\begin{equation} \label{eq:sinintS}
\lim_{\epsilon \to 0} \lim_{n \to \infty} \frac{\Sigma_1'}{16^n \log{n}/\sqrt{n}} = \frac{2\alpha}{\pi^{5/2}} \int_0^1 \frac{1}{\sqrt{1+\sinh^2(2\alpha t)}} \dif t.
\end{equation}
\end{lemma}
\begin{proof}
We show that the sums \eqref{eq:bigsumSabove} and \eqref{eq:bigsumSbelow} are asymptotically equal. This implies that $\frac{1}{2}\Sigma_1'$ and \eqref{eq:bigsumSabove} are asymptotically equal. We start with the upper bound. The inner sum in \eqref{eq:bigsumSabove} is smaller than
\begin{equation*}
\int_1^{\infty} e^{-\frac{1}{16} (j-1)^2 \epsilon^2 (1+m'(i,\epsilon)^2)} \dif j = \frac{2 \sqrt{\pi}}{\epsilon \sqrt{1+m'(i,\epsilon)^2}}.
\end{equation*}
Plugging this into \eqref{eq:bigsumSabove}, moving out all constants from the sum but keeping all $\epsilon$'s in it shows that it remains to evaluate
\begin{equation*}
\sum_{i=1}^{N} \frac{\epsilon}{\sqrt{1+m'(i,\epsilon)^2}}.
\end{equation*}
Again, we employ an integral estimate (using that $m'(i,\epsilon)$ is increasing on the interval $[0, N]$) to bound the last sum from above by
\begin{equation*}
\int_0^{N} \frac{\epsilon}{\sqrt{1+m'(i,\epsilon)^2}} \dif i =
\int_0^1 \frac{1}{\sqrt{ 1 + \sinh^2(2(t-\epsilon)\alpha)}} \dif t
\end{equation*}
after the substitution $t = i\epsilon$. As $\epsilon$ tends to $0$, this becomes the integral shown in \eqref{eq:sinintS}.

For the lower bound, the inner sum in \eqref{eq:bigsumSbelow} is at least
\begin{equation*}
\int_1^{N^2} e^{-\frac{1}{16} j^2 \epsilon^2 (1+M(i,\epsilon)^2)} \dif j = \frac{2 \sqrt{\pi}}{\epsilon \sqrt{1 + M(i, \epsilon)^2}}  \left( \erf\left( \frac{\sqrt{1+M(i,\epsilon)^2}}{4\epsilon} \right) - \erf\left( \frac{\epsilon\sqrt{(1+M(i,\epsilon)^2)}}{4} \right) \right).
\end{equation*}
Since the error function is monotonously increasing, and $M(i,\epsilon)$ is monotonously increasing on $[1, N]$ as well, the term involving the error functions is at least
\begin{equation*}
\erf\left( \frac{1}{4\epsilon} \right) - \erf\left( \frac{\epsilon}{2\sqrt{2}} \right)
\end{equation*}
which tends to $1$ as $\epsilon$ tends to $0$. We are left with the sum
\begin{equation*}
\sum_{i=1}^{N} \frac{\epsilon}{\sqrt{1 + M(i, \epsilon)^2}}
\end{equation*}
which is bounded from below by
\begin{equation*}
\int_1^{N} \frac{\epsilon}{\sqrt{1 + M(i, \epsilon)^2}} \dif i =
\int_{\epsilon}^1 \frac{1}{\sqrt{ 1 + \sinh^2(2\alpha x)}} \dif t
\end{equation*}
where $t = i \epsilon$. In the limit $\epsilon \to 0$ this again becomes the integral on the right-hand side in \eqref{eq:sinintS}.
\end{proof}

\section{Proof of Theorem~\ref{thm:mainthm}}
\label{sec:mainthm}

We are ready to prove Theorem \ref{thm:mainthm}. We first prove part (a) and then part (b).

\begin{proof}[Proof of Theorem \ref{thm:mainthm}(a)]
Whereas $\Sigma_1$ and $\Sigma_2$ depend on $\epsilon$, the total sum $\abs{R_{8n}}$ does not. In particular,
\begin{equation*}
\lim_{n \to \infty} \frac{\abs{R_{8n}}}{16^n \log{n}/\sqrt{n}} = \lim_{\epsilon \to 0} \lim_{n \to \infty} \frac{\abs{R_{8n}}}{16^n \log{n}/\sqrt{n}}.
\end{equation*}
The last double limit can be split in several pieces using that $\abs{R_{8n}} = 2^{2n} \binom{2n}{n} + 2\Sigma_1 + 2\Sigma_2 + 2\Sigma_3$. In particular, \eqref{eq:S3}, Lemma \ref{lem:S2}, and Lemma \ref{lem:S1} show that
\begin{equation*}
\lim_{\epsilon \to 0} \lim_{n \to \infty} \frac{\abs{R_{8n}}}{16^n \log{n}/\sqrt{n}} = \lim_{\epsilon \to 0} \lim_{n \to \infty} \frac{2 \Sigma_1}{16^n \log{n}/\sqrt{n}} = \frac{1}{2\pi^{3/2}} \int_0^1 \frac{1}{\sqrt{1 + \sin^2(\pi t/2)}} \dif t.
\end{equation*}
To evaluate the integral, substitute $x = \sin^4(\pi t/2)$ so that $2\pi \dif t = x^{-3/4} (1-\sqrt{x})^{-1/2} \dif x$. Hence
\begin{equation} \label{eq:intvalnew}
\int_0^1 \frac{1}{\sqrt{1 + \sin^2(\pi t/2)}} \dif t = \frac{1}{2\pi} \int_0^{1} x^{-3/4} (1-x)^{-1/2} \dif x = \frac{1}{2\pi} B \left(\frac{1}{4},\frac{1}{2}\right) = \frac{\Gamma(\frac{1}{4})\Gamma(\frac{1}{2})}{2\pi \Gamma(\frac{3}{4})}
\end{equation}
where $B$ is the beta function, which satisfies $B(m,n) = \Gamma(m)\Gamma(n)/\Gamma(m+n)$. Legendre's duplication formula for the gamma function yields $\Gamma(1/2) = \Gamma(1/4)\Gamma(3/4)/\sqrt{2\pi}$, showing that \eqref{eq:intvalnew} equals $\Gamma(\frac{1}{4})^2/\sqrt{8 \pi^3}$. This gives the desired result.

The skew-reciprocal case is entirely similar. With the same steps, we deduce
\begin{equation*}
\lim_{n \to \infty} \frac{\abs{S_{8n}}}{16^n \log{n}/\sqrt{n}} = \lim_{\epsilon \to 0} \lim_{n \to \infty} \frac{2 \Sigma'_1}{16^n \log{n}/\sqrt{n}} = \frac{4\alpha}{\pi^{5/2}} \int_0^1 \frac{1}{\sqrt{1+\sinh^2(2\alpha t)}} \dif t
\end{equation*}
where again $\alpha = \log{\sqrt{1+\sqrt{2}}}$. To evaluate the integral, substituting $x = \sinh(2\alpha t)$ yields $\dif x = 2\alpha\cosh(2\alpha t) \dif t = 2\alpha \sqrt{x^2+1} \dif t$. Therefore
\begin{equation*}
\int_0^1 \frac{1}{\sqrt{ 1 + \sinh^2(2\alpha t)}}  \dif t = \frac{1}{2\alpha} \int_0^1 \frac{1}{x^2+1} \dif x = \frac{1}{2\alpha} (\arctan(1) - \arctan(0)) = \frac{\pi}{8\alpha},
\end{equation*}
as claimed.
\end{proof}

\begin{proof}[Proof of Theorem \ref{thm:mainthm}(b)]
We prove the result for the reciprocal polynomials; an analogous argument works for the skew-reciprocals as well.
Write $n_0 = \frac{1}{4}(kr^2+ks^2+(-1)^{\frac{k+1}{2}})$. The second binomial coefficient in \eqref{eq:R8n-2} equals 
\begin{gather*}
\binom{2n-1}{n+n_0} = \left(\frac{1}{2} - \frac{n_0}{2n}\right) \binom{2n}{n+n_0} \quad \text{ if } k \equiv 1 \bmod{4}, \\
\binom{2n-1}{n+n_0-1} = \left( \frac{1}{2} + \frac{n_0}{2n} \right) \binom{2n}{n+n_0} \quad \text{ if } k \equiv 3 \bmod{4};
\end{gather*}
these identities also hold when $n_0 = n$. Therefore
\begin{equation*}
\abs{R_{8n-2}} = \frac{1}{2} \abs{R_{8n}} + \Sigma,
\end{equation*}
where
\begin{equation*}
\Sigma = \frac{1}{2n}  \sum_{(k,r,s) \in B_{5n}} \frac{1}{4} (1 + (-1)^{\frac{k+1}{2}}k(r^2+s^2)) \binom{2n}{n+\frac{1}{2}krs} \binom{2n}{n+n_0}.
\end{equation*}
Write
\begin{equation*}
V_t = \frac{1}{2n} \sum_{(k,r,s) \in B_{t}} \frac{1}{4} (1 + k(r^2+s^2)) \binom{2n}{n+\frac{1}{2}krs} \binom{2n}{n+n_0}.
\end{equation*}
Then $V_{5n}$ is at least as big as $\Sigma$ in absolute value. Estimating in each case the term
\begin{equation*}
\frac{1}{4} (1 + k(r^2+s^2))
\end{equation*}
by the maximum value it can possibly attain, we see that $V_{5n} - V_{m}$ is asymptotically at most $\Sigma_3$, whereas $V_m$ is asymptotically at most $\sqrt{\log{n}/n} (\Sigma_1 + \Sigma_2)$ (both up to a multiplicative constant). Both of these are negligible compared to $\abs{R_{8n}}$.
\end{proof}

\section{Proof and discussion of Theorem~\ref{thm:pyth}}
\label{sec:pyth}

In this section, we prove and discuss Theorem~\ref{thm:pyth}, starting with the proof.
\begin{proof}[Proof of Theorem~\ref{thm:pyth}]
Note that $X$ and $Y$ are the `randomised versions' of the numbers $b$ and $c$ defined in \eqref{eq:b} and \eqref{eq:c}. Hence the probability that $Y^2-X^2$ is a square equals the proportion of the total number of choices of the $a_i$ in the definitions of $c$ and $b$ that make $c^2-b^2$ a square. But this is precisely what is being counted in the proof of Proposition~\ref{prop:nrofsqdisc}, which is the number $\abs{R_{8n}}$ of reciprocal Littlewood polynomials of degree $8n$ with square discriminant; the only difference is that we were counting \emph{monic} Littlewood polynomials, meaning that the coefficient $a_{4n}$ was fixed to be $1$, whereas in the randomised case $a_{4n}$, or rather its counterpart $A_{4n}$, can be $-1$ as well. This only leads to multiplication by $2$ of the result obtained in Proposition~\ref{prop:nrofsqdisc}, because $c^2-b^2$ is a square if and only if $(-c)^2-b^2$ is a square -- that is, as many Pythagorean triples arise from \eqref{eq:b} and \eqref{eq:c} with $a_{4n} = 1$ as with $a_{4n} = -1$. However, in terms of \emph{proportions} this multiplication by $2$ does not matter, as there are twice as many tuples $(a_0, a_1, \ldots, a_{4n})$ when $a_{4n} = -1$ is allowed as well. All in all, we have  
\begin{equation*}
\mr{Prob}(Y^2-X^2=\square) = \frac{2\abs{R_{8n}}}{2^{4n+1}} = \frac{\abs{R_{8n}}}{2^{4n}}
\end{equation*}
and the result follows from Theorem~\ref{thm:mainthm}. The result for the probability of $Y^2+X^2$ being a square follows similarly from the case of skew-reciprocal Littlewood polynomials.
\end{proof}

We can think of Theorem~\ref{thm:pyth} as an asymptotic, Gaussian-weighted count of the number of Pythagorean triples with bounded legs or hypotenuse, because for any $\epsilon \in (0, 1/6)$ we have
\begin{equation} \label{eq:Gauss}
\mr{Prob}(Y^2-X^2 = \square) \sim \frac{1}{2\pi n} \sum_{\substack{\abs{\ell}, \abs{m} \leq n^{1/2+\epsilon} \\ \upsilon \in \{\pm 1\}}} e^{-\frac{\ell^2+m^2}{n}} \boldsymbol{1}_{(4m+\upsilon)^2-(4\ell)^2=\square}
\end{equation}
and similarly for $Y^2+X^2$. Indeed, observe that $\mr{Prob}(X=4\ell)=4^{-n}\binom{2n}{n+\ell}$ and $\mr{Prob}(Y=4m+\upsilon) = 4^{-n}\binom{2n}{n+m}/2$ for $\upsilon \in \{\pm 1\}$. Hence
\begin{align*}
\mr{Prob}(Y^2-X^2 = \square) &= \sum_{\substack{\abs{\ell}, \abs{m} \leq n \\ \upsilon \in \{\pm 1\}}} \mr{Prob}(X = 4\ell)\mr{Prob}(Y=4m + \upsilon)\boldsymbol{1}_{(4m+\upsilon)^2-(4\ell)^2=\square} \\
&= \frac{1}{2\cdot16^n} \sum_{\substack{\abs{\ell}, \abs{m} \leq n \\ \upsilon \in \{\pm 1\}}} \binom{2n}{n+\ell}\binom{2n}{n+m}\boldsymbol{1}_{(4m + \upsilon)^2-(4\ell)^2=\square} \\
&= \frac{1}{2\cdot16^n} \sum_{\substack{\abs{\ell}, \abs{m} \leq n^{1/2+\epsilon} \\ \upsilon \in \{\pm 1\}}} \left( \binom{2n}{n+\ell}\binom{2n}{n+m}\boldsymbol{1}_{(4m + \upsilon)^2-(4\ell)^2=\square}\right) + O(ne^{-2n^{2\epsilon}})
\end{align*}
which is asymptotic to the right-hand side in \eqref{eq:Gauss} by Proposition~\ref{prop:almostcentral}. The error term in the last line comes from the sum
\begin{align*}
\frac{1}{2\cdot 16^n} \sum_{\substack{n^{1/2+\epsilon} \leq \abs{\ell}, \abs{m} \leq n \\ \upsilon \in \{\pm 1\}}} \binom{2n}{n+\ell}\binom{2n}{n+m}\boldsymbol{1}_{(4m + \upsilon)^2-(4\ell)^2=\square} \leq \frac{4n^2}{2\cdot 16^n} \binom{2n}{n+n^{1/2+\epsilon}}^2 
\end{align*}
after again applying Proposition~\ref{prop:almostcentral}.

\section{Square discriminants in other degrees}
\label{sec:otherdeg}

In this section, we discuss Littlewood polynomials with square discriminant in degree $n \not \equiv 0, 6 \bmod{8}$.
The following surprising result, attributed to Alexei Entin in \cite[\S 4]{BK}, shows that such polynomials do not even exist in even degree $n \equiv 2, 4 \bmod{8}$.

\begin{lemma}[Entin] \label{lem:Entin}
Let $n \equiv 2, 4 \bmod{8}$ be a positive integer. Then no Littlewood polynomial of degree $n$ has square discriminant.
\end{lemma}
\begin{proof}
Suppose that $n$ is even and $f \in \F_{n}$. Set $p_n(X) = (X^{n+1}-1)/(X-1)$ and note that $f$ and $p_n$ coincide modulo $2$. Since $X^{n+1}-1$ and its derivative are coprime modulo $2$, the polynomial $p_n$ is separable over $\FF_2$. Thus $p_n$ is separable over the $2$-adic field $\Q_2$ as well by Hensel's lemma. The splitting field of $p_n$ over $\Q_2$, which is the cyclotomic extension $\Q_2(\zeta)/\Q_2$ where $\zeta$ is a primitive $n+1$-th root of unity, is an unramified extension of $\Q_2$ because $2$ and $n+1$ are coprime, see \cite[Proposition II.7.12]{Neukirch}. Writing $G(f/K)$ for the Galois group of $f$ over a field $K$, this implies that $G(p_n/\Q_2)$ is isomorphic to $G(p_n/\FF_2) = G(f/\FF_2) \leq G(f/\Q)$. The discriminant of $p_n$ is a square in $\Z_2$ if and only if it is $1 \bmod{8}$. A resultant calculation shows that $\Delta(p_n) = (-1)^{\frac{n(n-1)}{2}} (n+1)^{n-1}$, which is congruent to $5 \bmod{8}$ if $n \equiv 2, 4 \bmod{8}$ (and congruent to $1 \bmod{8}$ otherwise). Therefore $f$ cannot have square discriminant over $\Q$.  
\end{proof}

In the case of odd-degree Littlewood polynomials, the situation is different. Call a degree-$n$ polynomial $f$ \emph{nearly reciprocal} if $f(X) = \pm X^n f(X^{-1})$ and \emph{nearly skew-reciprocal} if $f(X) = \pm X^n f(-X^{-1})$. We give some examples:
\begin{itemize}
    \item Littlewood polynomials with vanishing square discriminant exist in any odd degree. Indeed, the nearly reciprocal polynomial given by 
    \begin{equation*}
    \qquad (X^{n+1}-1)(X^n + X^{n-1} + \cdots + X + 1) = (X-1)(X^n + X^{n-1} + \cdots + X + 1)^2 \in \F_{2n+1}
    \end{equation*} 
    has a multiple factor and thus its discriminant vanishes. 
    \item An odd-degree Littlewood polynomial with vanishing square discriminant is not necessarily nearly (skew-)reciprocal, or the product of such. Indeed, the polynomial
    \begin{equation*}
    \qquad (X + 1)^2 (X^2 - X + 1) (X^7 - X^5 + X^4 - X^3 + X^2 + 1)
    \end{equation*}
    has vanishing discriminant, but the Galois group of its splitting field is $C_2 \times S_7$.
    \item A computer experiment shows that all Littlewood polynomials of odd degree $\leq 29$ with nonvanishing square discriminant have a cyclotomic factor; in fact, each such polynomial is divisible by $X+1$ or $X-1$. Does there exist an odd-degree Littlewood polynomial without cyclotomic factors that has square discriminant? (If not, this would imply for example that no irreducible Littlewood polynomial of odd degree $n$ has Galois group contained in $A_n$.)  

    A related question, raised by Peled, Sen and Zeitouni \cite[\S 7]{PSZ}, is whether Littlewood polynomials with a repeated non-cyclotomic factor exist. The answer is `yes', the polynomial of degree $195$ given in \cite[Example~9]{DJS} apparently being the first known instance. In response to a question on MathOverflow \cite{Taylor}, Taylor found the example
    \begin{align*}
    \qquad \qquad (X^{18} + X^{16} + 2X^{15} + 2X^{13} + X^{12} + 2X^{11} + 3X^{10} + 3X^8 + 2X^7 &+ X^6 + 2X^5 + 2X^3 + 1) \\ & \!\!\!\!\!\!\!\!\times (X^2 + 1)(X - 1)(X^3 + X^2 - 1)^2
    \end{align*}
    of degree $27$.
\end{itemize}

\bibliographystyle{amsplain}

\begin{thebibliography}{10}

\bibitem{BBM}
L. Bary-Soroker, O. Ben-Porath, and V. Matei, \emph{Probabilistic {G}alois {T}heory -- {T}he {S}quare {D}iscriminant {C}ase}, preprint arXiv:2207.12493, 15 pp., 2022.

\bibitem{BKK}
L. Bary-Soroker, D. Koukoulopoulos, and G. Kozma, \emph{Irreducibility of random polynomials: general measures}, Invent. Math. \textbf{233} (2023), 1041--1120.

\bibitem{BK}
L. Bary-Soroker and G. Kozma, \emph{Irreducible polynomials of bounded height}, Duke Math. J. \textbf{169} (2020), 579--598.

\bibitem{BVPyth}
M. Benito and J. L. Varona, \emph{Pythagorean triangles with legs less than {$n$}}, J. Comput. Appl. Math. \textbf{143} (2002), 117--126.

\bibitem{Bhargava}
M. Bhargava, \emph{Galois groups of random integer polynomials and Van der Waerden's Conjecture}, preprint arXiv:2111.06507, 33 pp., 2021.

\bibitem{Borst_et_al}
C. Borst, E. Boyd, C. Brekken, S. Solberg, M. M. Wood, and P. M. Wood, \emph{Irreducibility of random polynomials}, Exp. Math. \textbf{27} (2018), 498--506.

\bibitem{BV}
E. Breuillard and P. Varj\'{u}, \emph{Irreducibility of random polynomials of large degree}, Acta Math. \textbf{223} (2019), 195--249.

\bibitem{Chela}
R. Chela, \emph{Reducible polynomials}, J. London Math. Soc. \textbf{38} (1963), 183--188.

\bibitem{DJS}
P. Drungilas, J. Jankauskas, and J. \v{S}iurys, \emph{On {L}ittlewood and {N}ewman polynomial multiples of {B}orwein polynomials}, Math. Comp. \textbf{87} (2018), 1523--1541.

\bibitem{Dubickas}
A. Dubickas, \emph{Salem numbers as {M}ahler measures of nonreciprocal units}, Acta Arith. \textbf{176} (2016), 81--88.

\bibitem{Durrett}
R. Durrett, \emph{Probability---theory and examples}, Cambridge Series in Statistical and Probabilistic Mathematics \textbf{49}, Cambridge University Press, Cambridge, 2019.

\bibitem{Erdelyi}
T. Erd\'{e}lyi, \emph{Do {F}lat {S}kew-{R}eciprocal {L}ittlewood {P}olynomials {E}xist?}, Constr. Approx. \textbf{56} (2022), 537--554.

\bibitem{Konyagin}
S. V. Konyagin, \emph{On the number of irreducible polynomials with {$0,1$} coefficients}, Acta Arith. \textbf{88} (1999), 333--350.

\bibitem{Kratzel}
E. Kr\"{a}tzel, \emph{Lattice points}, Mathematics and its Applications (East European Series) \textbf{33}, Kluwer Academic Publishers Group, Dordrecht, 1988.

\bibitem{Littlewood}
J. E. Littlewood, \emph{On polynomials {$\sum^{n}\pm z^{m}$}, {$\sum^{n}e^{\alpha _{m}i}z^{m}$}, {$z=e^{\theta_{i}}$}}, J. London Math. Soc. \textbf{41} (1966), 367--276.

\bibitem{MV}
H. L. Montgomery and R. C. Vaughan, \emph{Multiplicative number theory. {I}. {C}lassical theory}, Cambridge Studies in Advanced Mathematics \textbf{97}, Cambridge University Press, Cambridge, 2006.

\bibitem{Neukirch}
J. Neukirch, \emph{Algebraic Number Theory}, Grundlehren der Mathematischen Wissenschaften \textbf{322}, Springer-Verlag, Berlin, 1999.

\bibitem{NR}
W. G. Nowak and W. Recknagel, \emph{The distribution of {P}ythagorean triples and a three-dimensional divisor problem}, Math. J. Okayama Univ. \textbf{31} (1989), 213--220.

\bibitem{Nymann}
J. E. Nymann, \emph{On the probability that {$k$} positive integers are relatively prime {II}}, J. Number Theory \textbf{7} (1975), 406--412.

\bibitem{Odlyzko}
A. Odlyzko, \emph{Search for ultraflat polynomials with plus and minus one coefficients}, in: Connections in Discrete Mathematics, Cambridge Univ. Press, Cambridge, 2018, 39–55.

\bibitem{OW}
S. O'Rourke and P. M. Wood, \emph{Low-degree factors of random polynomials}, J. Theoret. Probab. \textbf{32} (2019), 1076--1104.

\bibitem{PSZ}
R. Peled, A. Sen and O. Zeitouni, \emph{Double roots of random {L}ittlewood polynomials}, Israel J. Math. \textbf{213} (2016), 55--77.

\bibitem{Sierpinski}
W. Sierpiński, \emph{O sumowaniu szeregu {$\sum_{n>a}^{n\leq b} \tau(n)f(n)$}, gdzie {$\tau(n)$} oznacza liczbę rozkładów liczby $n$ na sumę kwadratów dwóch liczb całkowitych} [Polish;~\emph{On the summation of the series {$\sum_{n>a}^{n\leq b} \tau(n)f(n)$}, where {$\tau(n)$} denotes the number of ways to write $n$ as the sum of squares of two integers}], Prace Mat.-Fiz. \textbf{18} (1907), 1--59. French in: \emph{Oeuvres choisies, Tome I}, PWN---\'{E}ditions Scientifiques de Pologne, Warsaw, 1974, 109--154.

\bibitem{Spencer}
J. Spencer with L. Florescu, \emph{Asymptopia}, Student Mathematical Library \textbf{71}, American Mathematical Society, Providence, RI, 2014.

\bibitem{Stronina}
M. I. Stronina, \emph{Integral points on circular cones}, Izv. Vys\v{s}. U\v{c}ebn. Zaved. Matematika \textbf{8} (1969), 112-116.

\bibitem{Taylor}
P. Taylor, \emph{Answer to question ``{M}ultiple roots of polynomials with coefficients {$\pm 1$}''}. Question posted by user Taras Banakh on MathOverflow, \url{https://mathoverflow.net/questions/424408/}, 2022.

\bibitem{VV}
P. Viana and P. M. Veloso, \emph{Galois theory of reciprocal polynomials}, Amer. Math. Monthly \textbf{109} (2002), 466-471.

\bibitem{VdW1}
B. L. van der Waerden, \emph{Die {S}eltenheit der {G}leichungen mit {A}ffekt}, Math. Ann. \textbf{109} (1934), 13--16.

\bibitem{VdW2}
B. L. van der Waerden, \emph{Die {S}eltenheit der reduziblen {G}leichungen und der {G}leichungen mit {A}ffekt}, Monatsh. Math. Phys. \textbf{43} (1936), 133--147.

\end{thebibliography}

\end{document}